\documentclass[10pt,twoside]{article}
\usepackage{amscd,amsfonts,amsmath,amssymb,amstext,amsthm}
\usepackage{mathpazo}
\usepackage{hyperref}
\usepackage{mathrsfs}

\def\RR{{\Bbb R}}
\def\CC{{\Bbb C}}
\def\BB{{\Bbb B}}

\def\PP{{\Bbb P}}
\def\HH{{\Bbb H}}

\def\ra{\rightarrow}

\newtheorem{theorem}{Theorem}[section]
\newtheorem{lemma}[theorem]{Lemma}
\newtheorem{proposition}[theorem]{Proposition}

\theoremstyle{remark}
\newtheorem{remark}[theorem]{Remark}

\newtheorem{example}[theorem]{Example}

\newcommand{\hq } {{/\kern -.185em/}}
\newcommand{\he } {{\kern -.050em\sim _H \kern -.050em }}
\newcommand{\eq } {{\kern -.100em\sim \kern -.100em}}
\newcommand{\eqs } {{\kern -.100em\sim }}
\date{}
\begin{document}
\title{Classical Symmetries of Complex Manifolds}

\author{Alan Huckleberry \& Alexander Isaev}

\maketitle

\sloppy
\parindent0ex
\parskip1.2ex
\pagenumbering {arabic}

\thispagestyle{empty}

\pagestyle{myheadings}
\markboth{Alan Huckleberry and Alexander Isaev}{Classical Symmetries of Complex Manifolds}

\begin{abstract} We consider complex manifolds that admit actions by holomorphic transformations of classical simple real Lie groups and classify all such manifolds in a natural situation. Under our assumptions, which require the group at hand to be dimension-theoretically large with respect to the manifold on which it is acting, our classification result states that the
manifolds which arise are described precisely as invariant open subsets of certain complex flag manifolds associated to the complexified groups.
\end{abstract}


\section{Introduction}
\setcounter{equation}{0}

In this paper we consider complex manifolds that admit \lq\lq classical symmetries\rq\rq, that is, manifolds $X$ endowed with almost effective actions by holomorphic transformations of classical simple real Lie groups $G_0$. Here we say that $G_0$ is a classical simple group if it is connected and its Lie algebra ${\mathfrak g}_0$ is a classical simple real Lie algebra, i.e. one of the following mutually exclusive possibilities holds: (a) ${\mathfrak g}_0$ is a real form of a classical simple complex Lie algebra; (b) ${\mathfrak g}_0$ is classical simple complex Lie algebra $\tilde{\mathfrak g}_0$ regarded as a real algebra (in this case we write ${\mathfrak g}_0=\tilde {\mathfrak g}_0^{\RR}$). Throughout the paper we refer to real Lie algebras of these two kinds as algebras of types I and II, respectively. Let ${\mathfrak g}$ be the complexification of ${\mathfrak g}_0$. For ${\mathfrak g}_0$ of type I the algebra ${\mathfrak g}$ is a classical simple complex Lie algebra, and for ${\mathfrak g}_0$ of type II it is isomorphic to $\tilde {\mathfrak g}_0\oplus\tilde {\mathfrak g}_0$. We will be primarily interested in type I algebras. There are numerous examples of manifolds with classical symmetries in this case, perhaps the best-known ones being irreducible Hermitian symmetric spaces corresponding to classical Lie algebras (see \cite{Hel}, p. 354). More general examples are given by open $G_0$-invariant subsets of complex flag manifolds $G/P$, where $G$ is the universal complexification of $G_0$ and $P$ is a parabolic subgroup of $G$. Here $G_0$ acts on $G/P$ by means of a covering map $G_0\ra\hat G_0$, where $\hat G_0\subset G$ is a real form of $G$. For results on the orbit structure of actions of real forms on complex flag manifolds we refer the reader to \cite{FHW}, \cite{Wo1}, \cite{Z}, \cite{ACM} and references therein.

We are interested in obtaining, whenever possible, a complete explicit classification of manifolds with classical symmetries for any choice of ${\mathfrak g}_0$. One motivation for our study is the following \lq\lq ball characterization theorem\rq\rq\, obtained in \cite{I1} (cf.\cite{BKS}). 

\begin{theorem}\label{ballcharthm} Let $X$ be a connected complex manifold of dimension $n\ge 2$, and $\BB^n$ the unit ball in $\CC^n$. Assume that the groups $\hbox{Aut}_{\mathcal O}(X)$ and $\hbox{Aut}_{\mathcal O}(\BB^n)$ of holomorphic automorphisms of $X$ and $\BB^n$ are isomorphic as topological groups endowed with the compact-open topology. Then $X$ is biholomorphic to either $\BB^n$ or to the complement of its closure in $\PP^n$. In particular, if in addition $X$ is Stein or hyperbolic, then $X$ is biholomorphic to $\BB^n$.
\end{theorem}

In this paper we will see that Theorem \ref{ballcharthm} is a corollary of general classification results for manifolds with classical symmetries. In fact, the conclusion of the theorem holds true if $X$ admits an almost effective action of the classical simple group $PSU_{n,1}:=SU_{n,1}/\hbox{(center)}$, which is isomorphic to $\hbox{Aut}_{\mathcal O}(\BB^n)$. Note that for $G_0=PSU_{n,1}$ we have ${\mathfrak g}_0={\mathfrak{su}}_{n,1}$ and ${\mathfrak g}={\mathfrak{sl}}_{n+1}(\CC)$.

Let $n:=\hbox{dim}\,X$ (we assume that $n\ge 2$), and define $k({\mathfrak g}_0)$ for ${\mathfrak g}_0$ of type I (resp. type II) to be the dimension of the defining complex representation of ${\mathfrak g}$ (resp. $\tilde {\mathfrak g}_0$). Observe that if $k({\mathfrak g}_0)\le n+1$ there is no reasonable classification of manifolds $X$ for type I algebras (resp. type II algebras) if ${\mathfrak g}$ (resp. $\tilde {\mathfrak g}_0$) belongs to the B- or D-series of simple complex Lie algebras (the B-series consists of ${\mathfrak{so}}_k(\CC)$ with $k=2l+1$, $l\ge 2$, and the D-series of ${\mathfrak{so}}_k(\CC)$ with $k=2l$, $l\ge 4$). Indeed, any direct product ${\mathcal Q}_{k({\mathfrak g}_0)-2}\times Z$, where ${\mathcal Q}_{m}$ for any $m\ge 1$ denotes the $m$-dimensional projective quadric in $\PP^{m+1}$ and $Z$ is any complex manifold of dimension $n-k({\mathfrak g}_0)+2$, admits an almost effective action of the orthogonal group $SO_{k({\mathfrak g}_0)}(\CC)$, and hence actions of all its real forms and the real Lie group $SO_{k({\mathfrak g}_0)}(\CC)^{\RR}$. The action is defined as the standard action on the first factor and the trivial action on the second factor.

For a similar reason, if $k({\mathfrak g}_0)\le n$ one cannot explicitly classify manifolds $X$ for type I algebras (resp. type II algebras) if ${\mathfrak g}$ (resp. $\tilde {\mathfrak g}_0$) belongs to the A- or C-series (the A-series consists of ${\mathfrak{sl}}_k(\CC)$ with $k\ge 2$, and the C-series of ${\mathfrak{sp}}_k(\CC)$ with $k=2l$, $l\ge 3$). Indeed, any direct product $\PP^{k({\mathfrak g}_0)-1}\times Z$, where $Z$ is an arbitrary complex manifold of dimension $n-k({\mathfrak g}_0)+1$, admits an almost effective action of the groups $SL_{k({\mathfrak g}_0)}(\CC)$, $Sp_{k({\mathfrak g}_0)}(\CC)$, again defined as the standard action on the first factor and the trivial action on the second factor. Even for $n=2$, ${\mathfrak g}_0={\mathfrak{sl}}_2(\RR)$ an explicit classification is unlikely to exist. Indeed, in \cite{I5} (see also \cite{I4}) we determined all hyperbolic 2-dimensional manifolds $X$, for which the group $\hbox{Aut}_{\mathcal O}(X)$ of holomorphic automorphisms of $X$ has dimension 3. The list of such manifolds is quite long, and one of the most difficult parts of the classification corresponds to the case when the Lie algebra of $\hbox{Aut}_{\mathcal O}(X)$ is isomorphic to ${\mathfrak{sl}}_2(\RR)$ (cf. \cite{I3}). It is possible that the results of \cite{I5} can be generalized to the case when a group $G_0$ with Lie algebra ${\mathfrak{sl}}_2(\RR)$ acts on $X$ almost effectively and properly, but it is unlikely that there exists an explicit classification if the assumption of properness is dropped. On the positive side we note, however, that in \cite{IK2} all manifolds that admit an {\it effective}\, action of $SU_n$ were determined (cf. \cite{U1}), and in \cite{I2} we showed that any hyperbolic manifold $X$ of dimension $n\ge 3$ with $\hbox{dim}\,\hbox{Aut}_{\mathcal O}(X)=\hbox{dim}\,{\mathfrak{sl}}_n(\CC)=n^2-1$ is biholomorphic to a product $\BB^{n-1}\times S$, where $S$ is a Riemann surface. 

Following the above discussion, in this paper we make the following assumption. We say that $k({\mathfrak g}_0)$ and $n$ satisfy Condition $(>)$ for ${\mathfrak g}_0$ of type I (resp. type II), if $k({\mathfrak g}_0)>n$ when ${\mathfrak g}$ (resp. $\tilde{\mathfrak g}_0$) belongs to the A- or C-series, and $k({\mathfrak g}_0)>n+1$ when ${\mathfrak g}$ (resp. $\tilde{\mathfrak g}_0$) belongs to the B- or D-series. Note that the algebra ${\mathfrak{sp}}_4(\CC)$ is isomorphic to ${\mathfrak{so}}_5(\CC)$ and we include it in the B-series, and the algebra ${\mathfrak{so}}_6(\CC)$ is isomorphic to ${\mathfrak{sl}}_4(\CC)$ and we include it in the A-series.  

It turns out that  Condition $(>)$ implies that $k({\mathfrak g}_0)=n+1$ for the case of the A- and C-series and $k({\mathfrak g}_0)=n+2$ for the case of the B- and D-series. We call these relations Condition $(=)$. This condition implies, in particular, that the maximal dimension of a classical simple complex group that can act on $X$ non-trivially is a polynomial in $n$. In contrast, it was shown in \cite{SW} that for non-semisimple complex groups the dimension can depend exponentially on $n$ even if $X$ is compact and homogeneous (but not K\"ahler).

We will now describe the content of the paper starting with algebras of type I. Our study of $G_0$-orbits in $X$ relies on passing to the universal complexification $G$ of $G_0$ (for the construction of $G$ see \cite{Ho}). The group $G$ acts on $X$ locally, but this action may not be globalizable. We overcome this difficulty by means of fibering every $G_0$-orbit $Y$ in $X$ over a $G_0$-orbit $\hat Y$ in a complex projective space $\PP^{\ell}$. Such a fibering $\varphi_{{\mathfrak g}_0,\hbox{\tiny $Y$}}:Y\ra \hat Y$ comes from the ${\mathfrak g}_0$-anticanonical fibration associated to the orbit $Y$ (we outline the construction and main properties of this fibration in Section \ref{ganticanonicalfib}). Since $\hbox{Aut}_{\mathcal O}(\PP^{\ell})$ is a complex Lie group, the universality property of $G$ implies that $G$ acts on $\PP^{\ell}$ globally and $\hat Y$ lies in a $G$-orbit $Z=G/H$. Moreover, the map $\varphi_{{\mathfrak g}_0,\hbox{\tiny $Y$}}$ extends to a holomorphic map $\psi_{\hbox{\tiny $S$}}$ defined on the union $S$ of local $G$-orbits of points in $Y$, and the set $\psi_{\hbox{\tiny $S$}}(S)$ is open in $Z$ (see Section \ref{ganticanonicalfib}). Therefore, $Z$ is at most $n$-dimensional, that is, $H$ has codimension at most $n$ in $G$ (see Remark \ref{codimensionofh}). It then follows that $H$ is a proper maximal parabolic subgroup of maximal dimension in $G$ and thus $Z$ is a complex flag manifold, on which $G_0$ acts as a real form of $G$. The parabolicity of $H$ is a consequence of a general result on dimension-theoretically maximal complex closed subgroups of complex Lie groups that we obtain in Proposition \ref{parabolic} in Section \ref{dimtheormax}. The proof of Proposition \ref{parabolic} utilizes Tits' normalizer fibration for complex groups. Note that for a complex group the corresponding anticanonical fibration coincides with the normalizer fibration (see \cite{HO}, p. 65).   

In fact, $Z$ is biholomorphic to $\PP^n$ if ${\mathfrak g}$ belongs to the A- or C-series, and to the $n$-dimensional projective quadric ${\mathcal Q}_n$ if ${\mathfrak g}$ belongs to the B- or D-series, except for the case of ${\mathfrak{so}}_5(\CC)$, where both $\PP^3$ and ${\mathcal Q}_3$ occur. Under the equivalence $Z\simeq\PP^n$ the group $G$ acts on $\PP^n$ as either $PSL_{n+1}(\CC)$ or $PSp_{n+1}(\CC)$ (embedded in $PSL_{n+1}(\CC)$ in the standard way), and under the equivalence $Z\simeq{\mathcal Q}_n$ it acts on ${\mathcal Q}_n$ as $PSO_{n+2}(\CC)$. Further, the map $\psi_{\hbox{\tiny $S$}}$ extends to a locally biholomorphic $G_0$-equivariant surjective map $\psi: X\ra U$, where $U$ is an open $G_0$-invariant subset of $Z$. In Sections \ref{aseries}--\ref{bseries} we go over all classical simple complex Lie algebras ${\mathfrak g}$ and their real forms ${\mathfrak g}_0$ and show that on all the occasions the map $\psi$ is in fact a biholomorphism. Thus, our classification of manifolds $X$ for type I algebras ${\mathfrak g}_0$ consists of all open connected $G_0$-invariant subsets of either $\PP^n$ or ${\mathcal Q}_n$. The determination of all open $G_0$-invariant subsets is an interesting question in its own right, especially for the actions of the real forms of $PSO_8(\CC)$ on ${\mathcal Q}_6$, where the triality property plays a role. To establish that $\psi$ is biholomorphic we analyze the specific orbit structure for the $G_0$-action on $Z$ in each case, so our proof of the biholomorphicity of $\psi$ is a case-by-case argument interwoven with the determination of open $G_0$-invariant subsets of $Z$. 

We now summarize the content of the preceding three paragraphs in the following theorem, which is the main result of the paper.

\begin{theorem}\label{maintheoremI} Let $X$ be a connected complex manifold of dimension $n\ge 2$ admitting an almost effective action by holomorphic transformations of a connected simple Lie group $G_0$. Let ${\mathfrak g}_0$ be the Lie algebra of $G_0$, and assume that ${\mathfrak g}_0$ is of type I. Suppose, furthermore, that $k({\mathfrak g}_0)$ and $n$ satisfy Condition $(>)$. Then $k({\mathfrak g}_0)$ and $n$ satisfy Condition (=), and $X$ is biholomorphic to a connected $G_0$-invariant open subset of the flag manifold $Z=G/H$, where $G$ is the universal complexification of $G_0$, the subgroup $H$ is a maximal parabolic subgroup of maximal dimension in $G$, and $G_0$ acts on $Z$ as a real form of $G$.
\end{theorem}      

Theorems \ref{supq}, \ref{slr}, \ref{slh} give a detailed classification of manifolds $X$ for the case of the A-series, Theorems \ref{sppq}, \ref{spr} for the case of the C-series, Theorems \ref{sopq}, \ref{sopqn=3}, \ref{sostar}, \ref{so62}, \ref{so71}, \ref{so53} for the case of the B- and D-series.

The above discussion for type I algebras holds for type II algebras as well, with the difference that $G_0$ acts on $Z$ as either $PSL_{n+1}(\CC)^{\RR}$ or $PSp_{n+1}(\CC)^{\RR}$ if $Z$ is biholomorphic to $\PP^n$, and as $PSO_{n+2}(\CC)^{\RR}$ if $Z$ is biholomorphic to ${\mathcal Q}_n$. All these actions are transitive on $Z$. This makes the case of type II algebras significantly easier than that of type I algebras and leads to the following result. 

\begin{theorem}\label{maintheoremII} Let $X$ be a connected complex manifold of dimension $n\ge 2$ admitting an almost effective action by holomorphic transformations of a connected simple Lie group $G_0$. Let ${\mathfrak g}_0$ be the Lie algebra of $G_0$, and assume that ${\mathfrak g}_0$ is of type II. Suppose, furthermore, that $k({\mathfrak g}_0)$ and $n$ satisfy Condition $(>)$. Then $k({\mathfrak g}_0)$ and $n$ satisfy Condition (=), and $X$ is biholomorphic to a flag manifold $Z=G/H$, where $G$ is the universal complexification of $G_0$, and $H$ is a maximal parabolic subgroup of maximal dimension in $G$.
\end{theorem}

Theorem \ref{maintheoremII} is proved in Section \ref{classsiftypeII}, where a detailed version of this theorem is also stated (see Theorem \ref{maintheoremIIdetails}).  

It should be stressed here that we utilize the ${\mathfrak g}_0$-anticanonical fibration in the proofs of Theorems \ref{maintheoremI}, \ref{maintheoremII} because the local action of $G$ on $X$ is not known to be globalizable a priori. If $G$ acted on $X$ globally, then using the anticanonical fibration would only be required in the proof of Proposition \ref{parabolic} for complex groups acting transitively. 

Before proceeding we note that smooth actions of various classical groups on special real manifolds of specific dimensions (that are often assumed to be odd) have been extensively studied. We refer the reader, for instance, to \cite{U1}--\cite{U3}, \cite{UK}, \cite{T} and references therein for details.

{\bf Acknowledgements.} This work was initiated while the first author was visiting the Australian National University. The research is supported by the Australian Research Council.

\vspace{-0.5cm}

\section{${\mathfrak g}_0$-Anticanonical Fibration}\label{ganticanonicalfib}
\setcounter{equation}{0}

In this section we introduce a tool that will be frequently used throughout the paper. More details on this subject can be found in \cite{HO}, \cite{Huck}.

Let $X$ be a connected complex manifold and $G_0$ a connected real Lie group acting on $X$ almost effectively by holomorphic transformations. Denote by $K^{-1}$ the anti-canonical line bundle over $X$, and let $\Gamma(X,K^{-1})$ be the complex vector space of holomorphic sections of $K^{-1}$. Clearly, $G_0$ acts on $\Gamma(X,K^{-1})$ by complex-linear transformations. Let $V$ be a complex $G_0$-stable finite-dimensional subspace of $\Gamma(X,K^{-1})$. Consider the meromorphic map
\begin{equation}
\varphi_{\hbox{\tiny $V$}}: X\ra\PP(V^*),\quad \varphi_{\hbox{\tiny $V$}}(x):=[f_x],\quad\hbox{with $f_x(\sigma):=\sigma(x)$}, \quad x\in X, \quad \sigma\in V,\label{mapvarphi}
\end{equation}
where $\PP(V^*)$ is the projectivization of the dual space $V^*$ of $V$, and $[a]\in\PP(V^*)$ denotes the equivalence class of $a\in V^*\setminus\{0\}$. The indeterminacy set of $\varphi_{\hbox{\tiny $V$}}$ coincides with the set of points $x$ for which $f_x\equiv 0$. Note that in the above formula the value $\sigma(x)$, if it is non-zero, is only well-defined for a particular choice of coordinates near $x$ and is multiplied by a constant independent of $\sigma$ when one passes to another coordinate system (this explains the need to consider $\PP(V^*)$ instead of $V^*$ in the definition of $\varphi_{\hbox{\tiny $V$}}$). Next, the action of $G_0$ on $V$ induces an action of $G_0$ on $V^*$ by complex-linear transformations, and hence an action on $\PP(V^*)$. The map $\varphi_{\hbox{\tiny $V$}}$ is $G_0$-equivariant with respect to this induced action.

One can also define $\varphi_{\hbox{\tiny $V$}}$ by assigning to every point $x\in X$ the set ${\mathcal H}_x:=\{\sigma\in V: \sigma(x)=0\}$, where the hypersurface ${\mathcal H}_x$ is regarded as an element of $\PP^(V^*)$.

Assume now that $G_0$ acts on $X$ transitively. In this case one can make a particular choice of $V$ for which $\varphi_{\hbox{\tiny $V$}}$ is holomorphic on all of $X$. Let $\Gamma_{\mathcal O}(X,TX)$ be the space of vector fields on $X$ which are holomorphic, i.e., ${\mathcal Z}\in \Gamma_{\mathcal O}(X,TX)$ if and only if
$$
{\mathcal Z}^{\CC}:=\displaystyle\frac{1}{2}\left({\mathcal Z}-iJ{\mathcal Z}\right)
$$
is a holomorphic $(1,0)$-field on $X$, where $J$ is the operator of almost complex structure on $X$. Denote by ${\mathfrak g}_0$ the Lie algebra of $G_0$ and for $v\in{\mathfrak g}_0$ let ${\mathcal X}_v\in\Gamma_{\mathcal O}(X,TX)$ be the corresponding complete holomorphic vector field on $X$, that is,
$$
\displaystyle {\mathcal X}_v(x):=\frac{d}{dt}\left[\exp(-tv)x\right]|_{t=0},\quad x\in X.
$$
Since the action of $G_0$ on $X$ is almost effective, the map $\iota: v\mapsto {\mathcal X}_v$ is a Lie algebra isomorphism onto its image. Next, consider the Lie subalgebra $\hat{\mathfrak g}_0:=\iota({\mathfrak g}_0)+J\iota({\mathfrak g}_0)\subset \Gamma_{\mathcal O}(X,TX)$. Clearly, $\hat{\mathfrak g}_0$ is a complex Lie algebra with the operator of complex structure induced by $J$. By the complexification ${\mathcal Z}\mapsto {\mathcal Z}^{\CC}$ we realize $\hat{\mathfrak g}_0$ as the algebra of holomorphic $(1,0)$-fields of the form ${\mathcal X}_{v_1}^{\CC}+i{\mathcal X}_{v_2}^{\CC}$ for $v_1,v_2\in{\mathfrak g}_0$.

Let $n:=\hbox{dim}\,X$. Set $V_{{\mathfrak g}_0}:=\bigwedge^n \hat{\mathfrak g}_0$. Clearly, $V_{{\mathfrak g}_0}$ is a finite-dimensional linear subspace of $\Gamma(X,K^{-1})$ and is spanned by sections of $\Gamma(X,K^{-1})$ of the form $\sigma={\mathcal Z}_1\wedge\dots\wedge {\mathcal Z}_n$, ${\mathcal Z}_j\in\hat{\mathfrak g}_0$. For every such $\sigma$ and every $g\in G_0$ the action of $g$ on $\sigma$ is given by $g\sigma=(g{\mathcal Z}_1)\wedge\dots\wedge(g{\mathcal Z}_n)$, where for ${\mathcal Z}={\mathcal X}_{v_1}^{\CC}+i{\mathcal X}_{v_2}^{\CC}$ we have
$$
g{\mathcal Z}:={\mathcal X}_{\hbox{\tiny Ad}g(v_1)}^{\CC}+i{\mathcal X}_{\hbox{\tiny Ad}g(v_2)}^{\CC}.
$$     
Since $g{\mathcal Z}$ lies in $\hat{\mathfrak g}_0$, it follows that $V_{{\mathfrak g}_0}$ is $G_0$-stable. Therefore, we can consider the corresponding map $\varphi_{{\mathfrak g}_0}:=\varphi_{\hbox{\tiny $V_{{\mathfrak g}_0}$}}$ (see (\ref{mapvarphi})). Since for every $x\in X$ there exists $\sigma={\mathcal Z}_1\wedge\dots\wedge {\mathcal Z}_n$ with $\sigma(x)\ne 0$, the map $\varphi_{{\mathfrak g}_0}$ is holomorphic on $X$. 

In coordinates the map $\varphi_{{\mathfrak g}_0}$ and its equivariance property can be described as follows. We let $m+1$ be the dimension of $V_{{\mathfrak g}_0}$, and fix a basis $\sigma_0,\dots,\sigma_m$ in $V_{{\mathfrak g}_0}$. Then for $x\in X$ we have $\varphi_{{\mathfrak g}_0}(x)=[\sigma_0(x):\dots:\sigma_m(x)]$ in homogeneous coordinates in $\PP^m$. Further, for $g\in G_0$ let $A_g\in GL_{m+1}(\CC)$ be the linear transformation by which $g$ acts on $V_{{\mathfrak g}_0}$, written in coordinates with respect to the basis $\sigma_0,\dots,\sigma_m$. Identifying $V_{{\mathfrak g}_0}$ and $\CC^{m+1}$ by means of these coordinates, we define an action of $G_0$ on $\CC^{m+1}$ as $g\mapsto (A_g^{-1})^T$. This action induces an action on $\PP^m$, and it is straightforward to verify that $\varphi_{{\mathfrak g}_0}(gx)=g\varphi_{{\mathfrak g}_0}(x)$ for all $x\in X$ and $g\in G_0$, i.e., $\varphi_{{\mathfrak g}_0}$ is $G_0$-equivariant. Therefore, $\varphi_{{\mathfrak g}_0}(X)$ is a $G_0$-orbit in $\PP^m$. Let $H_0$ and $J_0$ be the isotropy subgroups of a point $x_0\in X$ and the point $\varphi_{{\mathfrak g}_0}(x_0)\in\PP^m$, respectively. Clearly, $J_0$ contains $H_0$, and the corresponding fibration
$$
G_0/H_0\ra G_0/J_0,
$$
is called the ${\mathfrak g}_0$-anticanonical fibration. One can show that $J_0$ lies in the normalizer $N_{G_0}(H_0^{\circ})$ of the identity component $H_0^{\circ}$ of $H_0$ in $G_0$ (see \cite{HO}, pp. 61-62). An important property of the ${\mathfrak g}_0$-anticanonical fibration is that it is a holomorphic fiber bundle with fiber $(J_0/H_0^{\circ})/(H_0/H_0^{\circ})$, where $J_0/H_0^{\circ}$ is a complex Lie group (see \cite{HO}, p. 64).  

We now drop the assumption of the transitivity of the $G_0$-action on $X$. Let $Y\subset X$ be the $G_0$-orbit of a point $x_0\in X$. For every $y\in Y$ the tangent space $T_y(Y)$ to $Y$ at $y$ is spanned by the vectors ${\mathcal X}_v(y)$ with $v\in{\mathfrak g}_0$. Let $\hat T_y(Y)$ be the complex subspace in $T_y(X)$ spanned by the values of vector fields from $\hat{\mathfrak g}_0=\iota({\mathfrak g}_0)+J\iota({\mathfrak g}_0)$ at the point $y$. Since $G_0$ acts on $X$ by holomorphic transformations, the dimension of the maximal complex subspace of $T_y(Y)$ is independent of $y\in Y$. Hence the dimension of $\hat T_y(Y)$ is independent of $y$ as well, and we denote it by $\mu$.

Suppose that $\mu\ge 1$. Let ${\mathfrak g}:={\mathfrak g}_0\oplus i{\mathfrak g}_0$ be the complexification of ${\mathfrak g}_0$, and let ${\mathcal G}$ be any connected complex Lie group with Lie algebra ${\mathfrak g}$. The algebra ${\mathfrak g}$ acts on $X$ via the homomorphism $\tau: {\mathfrak g}\ra\hat{\mathfrak g}_0$, $\tau(v_1+iv_2):={\mathcal X}_{v_1}+J{\mathcal X}_{v_2}$. The map $\tau$ induces a local holomorphic action of ${\mathcal G}$ on $X$, and for $y\in Y$ we denote by $O_y$ the local ${\mathcal G}$-orbit of $y$. Clearly, $O_y$ is a complex $\mu$-dimensional submanifold of $X$, and $\hat T_y(Y)$ is the tangent space to $O_y$ at $y$. Let $S:=\cup_{y\in Y}O_y$. The set $S$ is a (possibly non-injectively) immersed complex submanifold of $X$ of dimension $\mu$ containing $Y$.   

As before, we realize $\hat{\mathfrak g}_0$ as the finite-dimensional complex Lie algebra of holomorphic $(1,0)$-fields on $X$ of the form ${\mathcal X}_{v_1}^{\CC}+i{\mathcal X}_{v_2}^{\CC}$ for $v_1,v_2\in{\mathfrak g}_0$. Set $V_{\hbox{\tiny $S$}}:=\bigwedge^{\mu} \hat{\mathfrak g}_0$. Clearly, $V_{\hbox{\tiny $S$}}$ is a finite-dimensional linear subspace of $\Gamma(X,K^{-1})$ and is spanned by sections of $\Gamma(X,K^{-1})$ of the form ${\mathcal Z}_1\wedge\dots\wedge {\mathcal Z}_{\mu}$, ${\mathcal Z}_j\in\hat{\mathfrak g}_0$. We let $\ell+1$ be the dimension of $V_{\hbox{\tiny $S$}}$ and fix a basis $\sigma_0,\dots,\sigma_\ell$ in $V_{\hbox{\tiny $S$}}$. Now for $x\in S$ we set
$$
\psi_{\hbox{\tiny $S$}}(x):=[\sigma_0(x):\dots:\sigma_{\ell}(x)]\in\PP^{\ell}.
$$
Clearly, $\psi_{\hbox{\tiny $S$}}$ is well-defined and holomorphic on $S$. As before, the $G_0$-action on $V_{\hbox{\tiny $S$}}$ induces an action of $G_0$ on $\CC^{\ell+1}$ by complex-linear transformations, and hence an action on $\PP^{\ell}$. It is straightforward to verify that $\psi_{\hbox{\tiny $S$}}(gy)=g\psi_{\hbox{\tiny $S$}}(y)$ for all $y\in Y$ and $g\in G_0$, that is, the map $\varphi_{{\mathfrak g}_0,\hbox{\tiny $Y$}}:=\psi_{\hbox{\tiny $S$}}|_Y$ is $G_0$-equivariant. Furthermore, the map $\psi_{\hbox{\tiny $S$}}$ is equivariant with respect to the local $G_0$-action on $S$. It follows that $\varphi_{{\mathfrak g}_0,\hbox{\tiny $Y$}}(Y)$ is a $G_0$-orbit in $\PP^{\ell}$. Let $H_0$ and $J_0$ be the isotropy subgroups of the points $x_0$ and $p_0:=\varphi_{{\mathfrak g}_0,\hbox{\tiny $Y$}}(x_0)\in\PP^{\ell}$, respectively. We have $H_0\subset J_0$, and the corresponding fibration
$$
G_0/H_0\ra G_0/J_0,
$$
is called the ${\mathfrak g}_0$-anticanonical fibration associated to the orbit $Y$. Arguing as in the proof of Proposition 1 on pp. 61-62 of \cite{HO}, we see that $J_0$ lies in $N_{G_0}(H_0^{\circ})$.

We now assume that the complexification ${\mathfrak g}$ of ${\mathfrak g}_0$ is semisimple, and let $G:=G_0^{\CC}$ be the universal complexification of $G_0$ (see \cite{Ho}). The group $G_0$ is mapped into $G$ by means of a homomorphism $\gamma$ such that $\gamma(G_0)$ is a closed real form of $G$. Since the center of any connected complex semisimple Lie group is finite, it follows that $\gamma$ is locally injective, hence the Lie algebra of $G$ is isomorphic to ${\mathfrak g}$. Further, as we have seen, the action of $G_0$ on $V_{\hbox{\tiny $S$}}$ induces an action of $G_0$ on $\PP^{\ell}$, that is, a homomorphism $\rho:G_0\ra\hbox{Aut}_{\mathcal O}(\PP^{\ell})\simeq PSL_{\ell+1}(\CC)$. Since $PSL_{\ell+1}(\CC)$ is a complex Lie group, the universality property of $G$ implies that there exists a complex Lie group homomorphism $\rho^{\CC}:G\ra PSL_{\ell+1}(\CC)$ such that
\begin{equation}
\rho=\rho^{\CC}\circ\gamma.\label{universal}
\end{equation}
In fact, since $\rho=\pi\circ\tilde\rho$, where $\tilde\rho: G_0\ra SL_{\ell+1}(\CC)$ is a linear representation and $\pi: SL_{\ell+1}(\CC)\ra PSL_{\ell+1}(\CC)$ is the natural factorization map, there exists a linear representation $\tilde\rho^{\CC}: G\ra SL_{\ell+1}(\CC)$ such that $\tilde\rho=\tilde\rho^{\CC}\circ\gamma$ and $\rho^{\CC}=\pi\circ\tilde\rho^{\CC}$. 

Thus, the group $G$ holomorphically acts on $\PP^{\ell}$ by way of $\rho^{\CC}$, and we have the inclusion of orbits $G_0.p_0\subset G.p_0$. Note that if $G$ is a simple group, either it acts on the orbit $G.p_0$ almost effectively, or $G.p_0=\{p_0\}$. In the latter case $J_0=G_0$, hence $H_0^{\circ}$ is normal in $G_0$. Since $G_0$ is simple, it follows that either $H_0=G_0$, or $H_0$ is discrete.   

Let $\gamma^{\CC}:{\mathcal U}\ra G$ be a continuation of $\gamma$ to a local complex Lie group homomorphism defined in a neighborhood ${\mathcal U}$ of $e\in {\mathcal G}$ such that $d_e\gamma^{\CC}(v_1+iv_2)=d_e\gamma(v_1)+{\mathcal J}d_e\gamma(v_2)$ for all $v_1,v_1\in{\mathfrak g}_0$, where $d_e$ denotes the differential at $e$ and ${\mathcal J}$ is the operator of complex structure in the Lie algebra of $G$. Clearly, $\gamma^{\CC}$ maps ${\mathcal U}$ onto a neighborhood of the identity in $G$. Then the local $G_0$-equivariance of $\psi_{\hbox{\tiny $S$}}$ and property (\ref{universal}) imply that for every $x\in S$ we have
\begin{equation}
\psi_{\hbox{\tiny $S$}}(gx)=\gamma^{\CC}(g)\psi_{\hbox{\tiny $S$}}(x),\label{univers1}
\end{equation} 
where $g\in {\mathcal G}$ is sufficiently close to $e$. In particular, $\psi_{\hbox{\tiny $S$}}(S)$ is an open subset of the orbit $G.p_0$. Since $\hbox{dim}\,S=\mu\le n$, it follows that $\hbox{dim}\,G.p_0\le n$. 

\begin{remark}\label{codimensionofh}Let $H$ be the isotropy subgroup of $p_0$ with respect to the $G$-action. Since $\hbox{dim}\,G.p_0\le n$, the codimension of $H$ in $G$ is at most $n$. This observation will be important for our future applications. In fact, in all situations considered below $H$ will turn out to be a proper dimension-theoretically maximal subgroup of $G$. Such subgroups are studied in the next section.
\end{remark}

Clearly, the above construction of the ${\mathfrak g}_0$-anticanonical fibration requires the dimension $\mu$ of $S$ to be positive and does not apply in the case when $x_0$ is a fixed point of the $G_0$-action on $X$. We will now state a simple lemma that will be useful for ruling out fixed points later in the paper. 

\begin{lemma}\label{fixedpoints} Let $G_0$ be a connected simple real Lie group acting almost effectively by holomorphic transformations on an $n$-dimensional complex manifold $X$. Assume that $G_0$ fixes a point in $X$ and has a positive-dimensional compact subgroup. Then the complexification ${\mathfrak g}$ of the Lie algebra ${\mathfrak g}_0$ of $G_0$ has a non-trivial complex $n$-dimensional representation.    
\end{lemma}

\begin{proof} Let $x_0$ be a fixed point of the $G_0$-action. Consider the isotropy representation of $G_0$ at $x_0$
$$
\alpha_{x_0}:G_0\ra GL(T_{x_0}(X),\CC),\quad g\mapsto d_{x_0}g,
$$
where $GL(T_{x_0}(X),\CC)\simeq GL_n(\CC)$ is the group of non-degenerate complex-linear transformations of the tangent space $T_{x_0}(X)$, and $d_{x_0}g$ denotes the differential at $x_0$ of the holomorphic automorphism by which an element $g\in G_0$ acts on $X$. Let $G$ be the universal complexification of $G_0$. Since $G$ is semisimple, its Lie algebra coincides with ${\mathfrak g}$. The universality property of $G$ implies that there exists a linear representation $\alpha_{x_0}^{\CC}:G\ra GL_n(\CC)$ such that $\alpha_{x_0}=\alpha_{x_0}^{\CC}\circ\gamma$, where $\gamma:G_0\ra G$ is a locally injective homomorphism such that $\gamma(G_0)$ is a closed real form of $G$.

Let $K_0$ be a maximal compact subgroup of $G_0$. The action of $K_0$ on $X$ is almost effective and can be linearized near $x_0$, which implies that the representation $\alpha|_{K_0}$ is almost faithful. Since $K_0$ is positive-dimensional and $G_0$ is simple, it follows that $\alpha$ is almost faithful, and therefore $\alpha_{x_0}^{\CC}|_{\gamma(G_0)}$ is almost faithful. Since $\gamma(G_0)$ is positive-dimensional, we see that $\hbox{ker}\,\alpha_{x_0}^{\CC}\ne G$, and thus ${\mathfrak g}$ has a non-trivial complex $n$-dimensional representation as required.   
\end{proof}  

\begin{remark}\label{fixedpointsrem} Observe that in Lemma \ref{fixedpoints} the group $G$ is simple if ${\mathfrak g}_0$ is of type I, and $G$ is a locally direct product $G=G_1\cdot G_2$, where each $G_j$ is a closed subgroup of $G$ with Lie algebra $\tilde{\mathfrak g}_0$, if ${\mathfrak g}_0$ is of type II. Hence for ${\mathfrak g}_0$ of type I the representation $\alpha_{x_0}^{\CC}$ is almost faithful and thus the induced representation of ${\mathfrak g}$ is faithful. If ${\mathfrak g}_0$ is of type II and the kernel of the representation of ${\mathfrak g}=\tilde{\mathfrak g}_0\oplus\tilde{\mathfrak g}_0$ is non-trivial, then this kernel coincides with one of the simple factors. We also remark that the statement of Lemma \ref{fixedpoints} holds true if $G_0$ is not necessarily semisimple, but is assumed instead to have no positive-dimensional normal closed subgroup $L$ such that the intersection $L\cap K_0$ is discrete. In this case $\alpha$ is almost faithful, and the proof given above applies.

\end{remark}

\section{Dimension-Theoretically Maximal Subgroups of Complex Lie Groups}\label{dimtheormax}
\setcounter{equation}{0}

Our arguments in the forthcoming sections rely on the fact that a proper dimension-theoretically maximal subgroup of a connected complex semisimple Lie group is parabolic. This result is obtained in the present section. In fact, we prove the following more general proposition, which -- at least in the semisimple case -- is well-known to specialists (cf. \cite{Wi}, p. 46, Lemma 1).     

\begin{proposition}\label{parabolic} Let $G$ be a complex connected Lie group and $H\subset G$ a proper closed complex subgroup. Assume that $H$ is dimension-theoretically maximal, that is, there exists no proper closed complex subgroup of $G$ of dimension greater than $\hbox{dim}\,H$. Then one of the following holds: 

(i) $G'\subset H$, and the Abelian group $G/H$ either has dimension 1, or is a compact torus without non-trivial proper subtori;

(ii) $G/H$ is biholomorphic to $\CC^p$ for $p\ge 1$;

(iii) $H$ contains the radical $R$ of $G$ and $H/R$ is a maximal parabolic subgroup of $G/R$ of maximal dimension. 

If, furthermore, $G$ is semisimple, then $H$ is a maximal parabolic subgroup of $G$ of maximal dimension.
\end{proposition}

\begin{proof} Let $N:=N_G(H^{\circ})$ be the normalizer of the identity component $H^{\circ}$ of $H$ in $G$. Since $H$ is dimension-theoretically maximal, we have either $N=G$, or $N^{\circ}=H^{\circ}$. We will consider these two cases separately.

{\bf Case 1.} Suppose first that $N=G$ and let $G_1:=G/H^{\circ}$. Since $H$ is dimension-theoretically maximal, $G_1$ has no positive-dimensional proper closed complex subgroups. Let $G_1=R_1\cdot S_1$ be the Levi-Malcev decomposition of $G$, where $R_1$ is the radical and $S_1$ is a semisimple Levi subgroup of $G_1$, respectively. Since $R_1$ is closed and connected, it either coincides with $G_1$ or is trivial. In the latter case $G_1$ is semisimple and hence contains a closed copy of $SL_2(\CC)$, which in turn contains a closed copy of $\CC^*$. This contradiction implies that $G_1=R_1$ is solvable. Since $G_1$ is solvable, there exists a positive-dimensional (not necessarily closed) complex Abelian subgroup $A$ in $G_1$. Let $C_{G_1}(A)$ be the centralizer of $A$ in $G_1$. Since $C_{G_1}(A)$ is a closed complex subgroup of $G_1$ containing $A$, we have $C_{G_1}(A)=G_1$. Hence $A$ is contained in the center $Z(G_1)$ of $G_1$, and therefore $Z(G_1)$ is positive-dimensional. Since $Z(G_1)$ is a closed complex subgroup of $G_1$, we have $Z(G_1)=G_1$, that is, $G_1$ is Abelian. Therefore, we obtain that $G'\subset H^{\circ}$, and hence $G_2:=G/H$ is a complex Abelian group that has no positive-dimensional proper closed complex subgroups.

By Theorem 3.2 of \cite{Morim}, any complex Abelian group is isomorphic to a direct product $Q\times\CC^m\times(\CC^*)^k$, where $Q$ is a complex Cousin group (i.e. $Q$ does not admit any non-constant holomorphic functions). It then follows that $G_2$ is either a Cousin group or 1-dimensional.

\begin{lemma}\label{cousin} A Cousin group that has no positive-dimensional proper closed complex subgroups is isomorphic to a compact torus.
\end{lemma}

\begin{proof} Let $Q$ be a Cousin group satisfying the assumptions of the lemma. We write it as $Q=\CC^{\ell}/\Gamma$ for some $\ell$, where $\Gamma$ is a non-trivial discrete subgroup in $\CC^{\ell}$. Let $V$ be the real subspace of $\CC^{\ell}$ spanned by $\Gamma$. Then $K:=V/\Gamma$ is the maximal compact subgroup of $Q$. Assuming that $Q$ is non-compact, we see that $V$ is a positive-dimensional proper real subspace of $\CC^{\ell}$. Since $Q$ has no positive-dimensional proper closed complex subgroups, $V$ is not a complex subspace. Let $W$ be a complement to $V\cap iV$ in $V$, such that the real Abelian group $W/W\cap\Gamma$ is compact (and hence is isomorphic to $(S^1)^s$ for some $s>0$). Then $(W+iW)/W\cap\Gamma$ is a complex closed subgroup of $Q$ isomorphic to $(\CC^*)^s$. Since $Q$ has no positive-dimensional proper closed complex subgroups, it follows that $Q$ is isomorphic to $(\CC^*)^s$, which contradicts the assumption that $Q$ is a Cousin group. Thus, we have shown that $Q$ is compact as required.      
\end{proof}

\vspace{-0.5cm}  
 
Lemma \ref{cousin} yields that the group $G_2$ either has dimension 1, or is a compact torus without non-trivial proper subtori. This is option (i) of the proposition.

{\bf Case 2.} Assume now that $N^{\circ}=H^{\circ}$ and consider the normalizer fibration $G/H\ra G/N$ with finite fiber $N/H$. This fibration coincides with the ${\mathfrak g}$-anticanonical fibration, where ${\mathfrak g}$ is the Lie algebra of $G$ (see \cite{HO}, p. 65). In particular, there is a homomorphism $\rho:G\ra PSL_{r+1}(\CC)=\hbox{Aut}_{\mathcal O}(\PP^r)$ such that $N=\rho^{-1}({\mathcal N})$, where ${\mathcal N}$ consists of all elements of $\rho(G)$ that fix some point $x_0$ in $\PP^r$. The homomorphism $\rho$ is the composition of a linear representation of $G$ in $GL_{r+1}(\CC)$ and the natural factorization map $GL_{r+1}(\CC)\ra PSL_{r+1}(\CC)$. By a result due to Chevalley \cite{C} (see also \cite{HO}, p. 31), the commutator subgroup $\rho(G)'=\rho(G')$ is an algebraic subgroup of $PSL_{r+1}(\CC)$. Hence the orbit $\rho(G)'.x_0$ is closed in $\rho(G).x_0$, and thus the subgroup $P:=\rho^{-1}(\rho(G)')N$ is closed in $G$. The dimension-theoretic maximality of $H$ now implies that either $P^{\circ}=N^{\circ}$, or $P=G$.

If $P^{\circ}=N^{\circ}$, then $G'\subset H^{\circ}$, which is again option (i) of the proposition. Let $P=G$. Then $\rho(G)=\rho(G)'{\mathcal N}$, that is, $\rho(G).x_0=\rho(G)'.x_0$. Let $\hat G$ be the Zariski closure of $\rho(G)$ in $PSL_{r+1}(\CC)$. Since $\hat G'=\rho(G)'$, it follows that the orbit $\rho(G)x_0$ is Zariski-closed in $\hat G.x_0$, and therefore $\rho(G).x_0=\hat G.x_0$. Let $\hat{\mathcal N}$ be the isotropy subgroup of $x_0$ under the $\hat G$-action. As a general remark, we now note that if $T$ is a closed connected normal subgroup of $\hat G$ such that $T.x_0$ is closed in $\hat G.x_0$, then $T\hat{\mathcal N}$ is a closed subgroup of $\hat G$, and thus $P_T:=\rho^{-1}(T\hat{\mathcal N}\cap\rho(G))$ is a closed subgroup of $G$ containing $N$. The dimension-theoretic maximality of $H$ now implies that either $P_T=G$, or $\hbox{dim}\,P_T=\hbox{dim}\,N$. If $P_T=G$, then $T.x_0=\rho(G).x_0=\hat G.x_0$, i.e., $T$ acts on $\rho(G).x_0$ transitively. If $\hbox{dim}\,P_T=\hbox{dim}\,N$, then $T$ acts on $\rho(G).x_0$ trivially, in particular $T\subset\hat{\mathcal N}$.   

We decompose $\hat G$ as $\hat G=R_u\rtimes (Z\cdot S)$, where $R_u$ is the unipotent radical, $S$ is a semisimple Levi subgroup, and $Z\simeq(\CC^*)^q$ for some $q\ge 0$. The radical $\hat R$ of $\hat G$ is then the subgroup $R_u\rtimes Z$. Since $\hat R$ is a connected normal algebraic subgroup of $\hat G$, the orbit $\hat R.x_0$ is closed in $\hat G.x_0$. Then we have either $P_{\hat R}=G$, or $\hbox{dim}\,P_{\hat R}=\hbox{dim}\,N$. We will now consider each of these two cases.  

{\bf Case 2.1.} Suppose first that $P_{\hat R}=G$. Since $R_u$ is a normal algebraic subgroup of $\hat G$, we have either $P_{R_u}=G$, or $\hbox{dim}\,P_{R_u}=\hbox{dim}\,N$. Assuming that $P_{R_u}=G$, we obtain $G/N=R_u/(R_u\cap\hat{\mathcal N})$. Since $R_u$ is a simply-connected nilpotent algebraic group acting on $\PP^r$ algebraically, it follows that $R_u\cap\hat{\mathcal N}$ has finitely many connected components and therefore is in fact connected. Thus $G/N\simeq\CC^p$ for some $p\ge1$, which is option (ii) of the proposition.

Assume now that $\hbox{dim}\,P_{R_u}=\hbox{dim}\,N$. In this case $R_u$ acts on $\rho(G).x_0$ trivially, thus $\hat G$ acts on this orbit as the group $Z\cdot S$. Since we assumed that $\hat R$ acts on $\rho(G).x_0$ transitively, it follows that $Z$ acts on the orbit transitively (in particular, $q>0$). This implies (e.g. for topological reasons) that $S$ acts on $\rho(G).x_0$ trivially, and therefore $G/N$ is biholomorphic to $(\CC^*)^s$ for some $s\ge 1$. Observe that each 1-dimensional factor of $Z$ acts on $\rho(G).x_0$ either transitively or trivially, hence there is a factor acting transitively, which implies that $G/N$ is in fact biholomorphic to $\CC^*$. Further, since both $R_u$ and $S$ lie in $\hat{\mathcal N}$ and the product $R_u\rtimes S$ contains the subgroup $\hat G'=\rho(G)'=\rho(G')$, we have $\rho(G')\subset{\mathcal N}$. It then follows that $G'$ is contained in $H^{\circ}$, and we are once again led to option (i) of the proposition.     

{\bf Case 2.2.} Suppose now that $\hbox{dim}\,P_{\hat R}=\hbox{dim}\,N$. In this case $\hat R\subset\hat{\mathcal N}$. Let $R$ be the radical of $G$. Clearly, $\rho(R)\subset\hat R$, and therefore $R\subset H^{\circ}$. Instead of the triple of groups $G$, $H$, $N$, one can now consider the triple $G/R$, $H/R$, $N/R$, and therefore from now on we assume that $G$ is semisimple. We will show that in this case $H$ is a maximal parabolic subgroup of maximal dimension in $G$, and thus obtain option (iii) of the proposition. 

Since all linear representations of a complex reductive group are algebraic, the homomorphism $\rho: G\ra PSL_{r+1}(\CC)$ arising from the normalizer fibration $G/H\ra G/N$ is algebraic. Suppose first that the orbit $\rho(G).x_0$ is closed in $\PP^r$. Then $\rho(G).x_0$ is a projective variety, and therefore $N$ is a parabolic subgroup of $G$ (in particular, $N$ is connected). Hence $H=N$ is a maximal parabolic subgroup of $G$ of maximal dimension.

Suppose now that $\rho(G).x_0$ is not closed in $\PP^r$, and let $Y_0$ be its closure. Since $\rho(G)$ is an algebraic subgroup of $PSL_{r+1}(\CC)$, it follows that $\rho(G).x_0$ is Zariski open in $Y_0$. In particular, the non-empty Zariski closed set $E_0:=Y_0\setminus\rho(G).x_0$ consists of $\rho(G)$-orbits of smaller dimensions. The dimension-theoretical maximality of $H$ implies that $\rho(G)$ fixes every point in $E_0$. 

Fix $\hat x_0\in E_0$ and consider the complex line $L_0\subset\CC^{r+1}$ defined by $\hat x_0$. Recall that $\rho=\pi\circ\tilde\rho$, where $\tilde\rho:G\ra GL_{r+1}(\CC)$ is a linear representation of $G$, and $\pi:GL_{r+1}(\CC)\ra PSL_{r+1}(\CC)$ is the natural projection. Clearly, the line $L_0$ is $\tilde\rho(G)$-invariant. Let $L_0':=\CC^{r+1}/L_0$ and define $\sigma_0:\PP^r\setminus\{\hat x_0\}\ra\PP L_0'$ to be the projectivization of the projection $\CC^{r+1}\setminus L_0\ra L_0'\setminus\{0\}$. Observe that for every $y\in\PP L_0'$ the set $\sigma_0^{-1}(y)\cup\{\hat x_0\}$ is a projective line in $\PP^r$. The map $\sigma_0$ is $\rho(G)$-equivariant, hence $\sigma_0(\rho(G).x_0)=\rho(G).x_1$ for some $x_1\in\PP L_0'$. The dimension-theoretical maximality of $H$ now implies that either $\hbox{dim}\,\rho(G).x_1=\hbox{dim}\,\rho(G).x_0$, or $\rho(G).x_1=\{x_1\}$. The second case leads to a contradiction since it follows in this case that $\rho(G).x_0$ is a positive-dimensional non-compact orbit of a complex semisimple Lie group lying in a projective line. Therefore we have $\hbox{dim}\,\rho(G).x_1=\hbox{dim}\,\rho(G).x_0$. It then follows that the restriction of $\sigma_0$ to $\rho(G).x_0$ is a proper finite-to-one map. This implies that $\rho(G).x_1$ is not closed in $\PP L_0'$, and we can repeat the above argument for $\rho(G).x_1$ in place of $\rho(G).x_0$ and $\PP L_0'$ in place of $\PP^r$. Iterating this argument sufficiently many times we arrive at a situation when an orbit $\rho(G).x_k$ with $\hbox{dim}\,\rho(G).x_k=\hbox{dim}\,\rho(G).x_0$ lies in a projective subspace $\PP L_{k-1}'$ of dimension less than $\hbox{dim}\,\rho(G).x_0$. This contradiction finalizes the proof of the parabolicity of $H$.

The proof of the proposition is complete.
\end{proof}

We apply Proposition \ref{parabolic} mainly to simple complex groups $G$ whose Lie algebra ${\mathfrak g}$ is a classical matrix algebra. In this case the subgroup $H$ is a maximal parabolic subgroup of maximal dimension in $G$. Let us now recall the well-known description of such subgroups for each of the series A, B, C, D. Note that proper maximal parabolic subgroups of a simple complex Lie group correspond to omitting one node in its Dynkin diagram (see e.g. \cite{FH}, pp. 383--395). 

Let $\hat G$ be the classical matrix group (one of $SL_k(\CC)$, $Sp_k(\CC)$, $SO_k(\CC)$) with Lie algebra ${\mathfrak g}$. The group $G$ covers the quotient $G/\hbox{(center)}$, which is also covered by $\hat G$. Accordingly, the parabolic subgroups of $G$ are obtained from those of $\hat G$ by taking the direct and inverse images under the respective covering maps. Therefore, it is sufficient to describe proper maximal parabolic subgroups of the classical matrix groups.

{\bf A-series.} Let $\hat G=SL_k(\CC)$, $k\ge 3$. Every maximal parabolic subgroup of $\hat G$ is conjugate to the subgroup that stabilizes the $m$-dimensional subspace 
$$
L_{m,k}^1:=\{z_{m+1}=0,\dots,z_k=0\}\subset\CC^k,
$$
for some $1\le m\le k-1$. There are exactly two conjugacy classes of maximal parabolic subgroups of maximal dimension in $\hat G$; they correspond to $L_{1,k}^1$ and $L_{k-1,k}^1$. The quotient of $\hat G$ by any of these subgroups is biholomorphic to $\PP^{k-1}$. Let $P$ be a parabolic subgroup of $G$ arising from either conjugacy class in $\hat G$. Then $G/P$ is biholomorphically equivalent to $\PP^{k-1}$, and, since $P$ contains the center of $G$, under this equivalence the group $G$ acts on $\PP^{k-1}$ as the full automorphism group $\hbox{Aut}_{\mathcal O}(\PP^n)=PSL_k(\CC):=SL_k(\CC)/\hbox{(center)}$.

{\bf C-series.} Let $\hat G=Sp_k(\CC)$, $k=2l$, $l\ge 3$. We realize $\hat G$ as the group of all matrices in $GL_k(\CC)$ preserving the skew-symmetric bilinear form 
\begin{equation}
(z,w):=\sum_{j=1}^{l}\left(z_jw_{l+j}-z_{l+j}w_j\right),\label{matrixr}
\end{equation}
where $z:=(z_1,\dots,z_k)$, $w:=(w_1,\dots,w_k)$ are vectors in $\CC^k$. Every maximal parabolic subgroup of $\hat G$ is conjugate to the subgroup that stabilizes the $m$-dimensional subspace
$$
L_{m,k}^2:=\{z_{m+1}=0,\dots,z_k=0\}\subset\CC^k,
$$
for some $1\le m\le l$. There is only one conjugacy class of maximal parabolic subgroups of maximal dimension in $\hat G$; it corresponds to $L_{1,k}^2$. The quotient of $\hat G$ by any of these subgroups is biholomorphic to $\PP^{k-1}$. For any parabolic subgroup $P$ of $G$ arising from this conjugacy class in $\hat G$, the quotient $G/P$ is biholomorphically equivalent to $\PP^{k-1}$, and the equivalence can be chosen so that $G$ acts on $\PP^{k-1}$ as the group $PSp_k(\CC):=Sp_k(\CC)/\hbox{(center)}$ embedded into $PSL_k(\CC)$ in the standard way. 
   
{\bf B-series.} Let $\hat G=SO_k(\CC)$, $k=2l+1$, $l\ge 2$. We realize $\hat G$ as the group of all matrices $A\in SL_k(\CC)$ such that $A=(A^T)^{-1}$. Every maximal parabolic subgroup of $\hat G$ is conjugate to the subgroup that stabilizes the $m$-dimensional subspace
$$
L_{m,l}^3:=\{(z_1,iz_1,\dots,z_m,iz_m,0,\dots,0): z_1,\dots,z_m\in\CC\}\subset\CC^k,
$$
for some $1\le m\le l$. For $l\ge 3$ there is only one conjugacy class of maximal parabolic subgroups of maximal dimension in $\hat G$; it corresponds to $L_{1,l}^3$. The quotient of $\hat G$ by any of these subgroups is biholomorphic to the $(k-2)$-dimensional projective quadric 
$$
{\mathcal Q}_{k-2}:=\left\{\zeta\in\PP^{k-1}:\zeta_1^2+\dots+\zeta_k^2=0\right\},
$$
where $\zeta:=(\zeta_1:\dots:\zeta_k)$ are homogeneous coordinates in $\PP^{k-1}$. For any parabolic subgroup $P$ of $G$ arising from this conjugacy class in $\hat G$, the quotient $G/P$ is biholomorphically equivalent to ${\mathcal Q}_{k-2}$, and under this equivalence $G$ acts on ${\mathcal Q}_{k-2}$ as the full automorphism group $\hbox{Aut}_{\mathcal O}({\mathcal Q}_{k-2})=SO_k(\CC)$. In contrast, for $l=2$ there are exactly two conjugacy classes of maximal parabolic subgroups of maximal dimension in $\hat G$. The quotient of $\hat G$ by any subgroup corresponding to $L_{1,2}^3$ is biholomorphically equivalent to ${\mathcal Q}_3$, and for a parabolic subgroup $P$ of $G$ arising from this conjugacy class, under the equivalence $G/P\simeq{\mathcal Q}_3$ the group $G$ acts on ${\mathcal Q}_3$ as $SO_5(\CC)$. Further, the quotient of $\hat G$ by any subgroup corresponding to $L_{2,2}^3$ is biholomorphic to $\PP^3$, and for a parabolic subgroup $P$ of $G$ arising from this conjugacy class, the equivalence $G/P\simeq\PP^3$ can be chosen so that $G$ acts on $\PP^3$ as the group $PSp_4(\CC)$ embedded into $PSL_4(\CC)$ in the standard way. 

{\bf D-series.} Let $\hat G=SO_k(\CC)$, $k=2l$, $l\ge 4$. Again, we realize $\hat G$ as the group of all matrices $A\in SL_k(\CC)$ such that $A=(A^T)^{-1}$. Every maximal parabolic subgroup of $\hat G$ is conjugate to either the subgroup that stabilizes the $m$-dimensional subspace
$$
L_{m,l}^4:=\{(z_1,iz_1,\dots,z_m,iz_m,0,\dots,0): z_1,\dots,z_m\in\CC\}\subset\CC^k,
$$
for some $1\le m\le l$, or to the subgroup that stabilizes the $l$-dimensional subspace
$$
L_{l,l}^{4'}:=\{(-z_1,iz_1,z_2,iz_2,\dots,z_l,iz_l): z_1,\dots,z_l\in\CC\}\subset\CC^k.
$$
For $l\ge 5$ there is only one conjugacy class of maximal parabolic subgroups of maximal dimension in $\hat G$; it corresponds to $L_{1,l}^4$. The quotient of $\hat G$ by any of these subgroups is biholomorphic to ${\mathcal Q}_{k-2}$. For any parabolic subgroup $P$ of $G$ arising from this conjugacy class in $\hat G$, the quotient $G/P$ is biholomorphically equivalent to ${\mathcal Q}_{k-2}$, and under this equivalence $G$ acts on ${\mathcal Q}_{k-2}$ as the group $PSO_k(\CC):=SO_k(\CC)/\hbox{(center)}=\hbox{Aut}_{\mathcal O}({\mathcal Q}_{k-2})^{\circ}$. In contrast, for $l=4$ there are exactly three conjugacy classes of maximal parabolic subgroups of maximal dimension in $\hat G$; they correspond to $L_{1,4}^4$, $L_{4,4}^4$, $L_{4,4}^{4'}$. The quotient of $\hat G$ by any of these subgroups is biholomorphic to ${\mathcal Q}_6$. Let $P$ be a parabolic subgroup of $G$ arising from any of the three conjugacy classes in $\hat G$. Then $G/P$ is biholomorphically equivalent to ${\mathcal Q}_6$, and under this equivalence $G$ acts on ${\mathcal Q}_6$ as the group $PSO_8(\CC)$.

For more details on parabolic subgroups see, for example, \cite{Z}, pp. 93--95.

From now on we adopt the notation introduced in Section \ref{ganticanonicalfib}. In the next section we deal with the easier case of type II algebras and classify the corresponding manifolds $X$.

\section{Proof of Theorem \ref{maintheoremII}}\label{classsiftypeII}
\setcounter{equation}{0}

Since ${\mathfrak g}=\tilde{\mathfrak g}_0\oplus\tilde{\mathfrak g}_0$, the group $G$ is represented as a locally direct product $G=G_1\cdot G_2$, where $G_j$ is a complex closed normal subgroup of $G$ with Lie algebra $\tilde{\mathfrak g}_0$. Maximal proper parabolic subgroups of maximal dimension in $G$ are then described as $P_1\cdot G_2$ and $G_1\cdot P_2$, where $P_j$ is a maximal proper parabolic subgroup of maximal dimension in $G_j$, $j=1,2$.

Suppose first that $\tilde{\mathfrak g}_0$ is either ${\mathfrak{sl}}_k(\CC)$ or ${\mathfrak{sp}}_k(\CC)$, where in the latter case $k=2l$, $l\ge 3$. It follows from the description of parabolic subgroups for the A- and C-series given at the end of Section \ref{dimtheormax} that $P_j$ has codimension $k-1$ in $G$. Choose $Y$ to be any $G_0$-orbit in $X$. Since by Remark \ref{codimensionofh} we have $\hbox{codim}\,H\le n$, Proposition \ref{parabolic} yields that $k=n+1$ and $H$ is a maximal parabolic subgroup of maximal dimension in $G$ with $G/H$ biholomorphic to $\PP^n$. Under the equivalence $G/H\simeq\PP^n$ the group $G$ acts on $G/H$ as $PSL_{n+1}(\CC)$ in the case of the A-series and as $PSp_{n+1}(\CC)$ embedded in $PSL_{n+1}(\CC)$ in the standard way in the case of the C-series. It then follows that $G_0$ acts on $G/H\simeq\PP^n$ as the group $PSL_{n+1}(\CC)^{\RR}$ in the case of the A-series and as the group $PSp_{n+1}(\CC)^{\RR}$ in the case of the C-series. Since the action of each of $PSL_{n+1}(\CC)^{\RR}$ and $PSp_{n+1}(\CC)^{\RR}$ is transitive on $\PP^n$, the orbit $Y$ is open in $X$. Thus, $Y=X$ and $\varphi_{{\mathfrak g}_0,\hbox{\tiny $Y$}}:X\ra\PP^n$ is a covering map. Hence $X$ is biholomorphic to $\PP^n$.

Suppose next that $\tilde{\mathfrak g}_0={\mathfrak{so}}_k(\CC)$, with $k\ge 5$, $k\ne 6$. It follows from the description of parabolic subgroups for the B- and D-series given in Section \ref{dimtheormax} that $P_j$ has codimension $k-2$ in $G$. Choose $Y$ to be any $G_0$-orbit in $X$. Since $\hbox{codim}\,H\le n$, Proposition \ref{parabolic} yields that $k=n+2$ and $H$ is a maximal parabolic subgroup of maximal dimension in $G$. Hence $G/H$ is biholomorphic to ${\mathcal Q}_n$ if $k\ge 7$, and to either $\PP^3$ or ${\mathcal Q}_3$ if $k=5$. Under the equivalence $G/H\simeq{\mathcal Q}_n$ the group $G$ acts on $G/H$ as $PSO_{n+2}(\CC)$, and if $k=5$ and $G/H$ is equivalent to $\PP^3$, then under this equivalence $G$ acts on $G/H$ as $PSp_4(\CC)$ embedded in $PSL_4(\CC)$ in the standard way. It then follows that for $k\ge 7$ the group $G_0$ acts on $G/H\simeq{\mathcal Q}_n$ as the group $PSO_{n+2}(\CC)^{\RR}$, and for $k=5$ it acts on either $G/H\simeq{\mathcal Q}_3$ or $G/H\simeq\PP^3$ as either the group $SO_5(\CC)^{\RR}$ or the group $PSp_4(\CC)^{\RR}$, respectively. Since all these actions are transitive, $Y$ is open in $X$. Thus, $Y=X$ and $\varphi_{{\mathfrak g}_0,\hbox{\tiny $Y$}}$ is a covering map. The proof of Theorem \ref{maintheoremII} is complete.

We will now restate Theorem \ref{maintheoremII} in more detail as follows.

\begin{theorem}\label{maintheoremIIdetails} Let $X$ be a connected complex manifold of dimension $n\ge 2$ admitting an almost effective action by holomorphic transformations of a connected simple Lie group $G_0$. Let ${\mathfrak g}_0$ be the Lie algebra of $G_0$, and assume that ${\mathfrak g}_0$ is of type II. Suppose, furthermore, that $k({\mathfrak g}_0)$ and $n$ satisfy Condition $(>)$. Then $k({\mathfrak g}_0)$ and $n$ satisfy Condition (=), and the following holds:

(i) if $\tilde{\mathfrak g}_0$ is one of ${\mathfrak{sl}}_{n+1}(\CC)$, ${\mathfrak{sp}}_{n+1}(\CC)$, where in the latter case $n$ is odd, then $X$ is biholomorphic to $\PP^n$;

(ii) if $\tilde{\mathfrak g}_0={\mathfrak{so}}_{n+2}(\CC)$ for $n\ge 5$, then $X$ is biholomorphic to ${\mathcal Q}_n$; if $\tilde{\mathfrak g}_0={\mathfrak{so}}_5(\CC)$ (here $n=3$), then $X$ is biholomorphic to either $\PP^3$ or ${\mathcal Q}_3$.
\end{theorem}

\section{Preparation for the Proof of Theorem \ref{maintheoremI}}\label{preparation}
\setcounter{equation}{0}

In this section we begin the proof of Theorem \ref{maintheoremI} and obtain the following preliminary result.

\begin{proposition}\label{maintheoremIprep} Let $X$ be a connected complex manifold of dimension $n\ge 2$ admitting an almost effective action by holomorphic transformations of a connected simple Lie group $G_0$. Let ${\mathfrak g}_0$ be the Lie algebra of $G_0$, and assume that ${\mathfrak g}_0$ is of type I. Suppose, furthermore, that $k({\mathfrak g}_0)$ and $n$ satisfy Condition $(>)$. Then $k({\mathfrak g}_0)$ and $n$ satisfy Condition (=), and for every $G_0$-orbit $Y$ in $X$ the following holds:

(i) $H$ is maximal parabolic subgroup of maximal dimension in $G$;

(ii) $S$ is open and the map $\psi_{\hbox{\tiny $S$}}$ extends to a $G_0$-equivariant locally biholomorphic surjective map $\psi: X\ra U$, where $U$ is an open $G_0$-invariant subset of $G.p_0\simeq G/H$ and $G_0$ acts on $G/H$ as a real form of $G$.
\end{proposition}

\begin{proof} We consider two cases.

{\bf Case 1.} Assume first that ${\mathfrak g}$ is either ${\mathfrak{sl}}_k(\CC)$ or ${\mathfrak{sp}}_k(\CC)$, where in the latter case $k=2l$, $l\ge 3$. We start by observing that the $G_0$-action on $X$ is fixed point free. Indeed, for $n\ge 2$ the group $G_0$ contains a non-trivial compact subgroup. Hence if $G_0$ has a fixed point in $X$, by Lemma \ref{fixedpoints} and Remark \ref{fixedpointsrem} the algebra ${\mathfrak g}$ has a faithful $n$-dimensional representation. However, neither of ${\mathfrak{sl}}_k(\CC)$, ${\mathfrak{sp}}_k(\CC)$ has faithful representations in dimensions less than $k$. This follows, for example, from Weyl's dimension formula and the well-known dimensions of the fundamental representations of simple complex Lie algebras (see e.g. \cite{S}, Table 1). Thus, $G_0$ has no fixed points in $X$.

Let $Y$ be any $G_0$-orbit in $X$. Observe that $G$ acts on the corresponding orbit $G.p_0\subset\PP^{\ell}$ almost effectively, that is, $G.p_0$ is positive-dimensional. Indeed, if $G.p_0=\{p_0\}$, the subgroup $H_0$ is discrete, and therefore $\hbox{dim}\,Y=\hbox{dim}\,G_0=k^2-1$ if ${\mathfrak g}={\mathfrak{sl}}_k(\CC)$, and $\hbox{dim}\,Y=\hbox{dim}\,G_0=2l^2+l$ if ${\mathfrak g}={\mathfrak{sp}}_k(\CC)$. Neither of these identities can hold, however, since $\hbox{dim}\,Y\le n$. 

Further, recall that by Remark \ref{codimensionofh} we have $\hbox{codim}\,H\le n$. On the other hand, if $P$ is a maximal parabolic subgroup of $G$, then $\hbox{codim}\,P\ge k-1$ (see the end of Section \ref{dimtheormax}). By Proposition \ref{parabolic} we then obtain that $k=n+1$ and $H$ is a maximal parabolic subgroup of maximal dimension in $G$. It then follows that $G/H$ is biholomorphic to $\PP^n$. Under the equivalence $G/H\simeq\PP^n$ the group $G$ acts on $G/H$ as the group $PSL_{n+1}(\CC)$ if ${\mathfrak g}={\mathfrak{sl}}_{n+1}(\CC)$, and as the group $PSp_{n+1}(\CC)$ embedded in $PSL_{n+1}(\CC)$ in the standard way if ${\mathfrak g}={\mathfrak{sp}}_{n+1}(\CC)$.  

Next, since $\psi_{\hbox{\tiny $S$}}(S)$ is open in the $n$-dimensional manifold $G.p_0\simeq\PP^n$, it follows that $S$ is open as well. Since this is the case for every $G_0$-orbit in $X$, the map $\psi_{\hbox{\tiny $S$}}$ extends to a $G_0$-equivariant holomorphic map $\psi: X\ra\PP^{\ell}$. Clearly, for any $G_0$-orbit $Y'$ in $X$ and the corresponding set $S'$ associated to $Y'$, we have $\psi|_{S'}=\psi_{\hbox{\tiny $S'$}}$, that is, $\psi$ is a continuation to all of $X$ of every map $\varphi_{{\mathfrak g}_0,\hbox{\tiny $Y'$}}$. Therefore, (\ref{univers1}) implies
\begin{equation}
\psi(gx)=\gamma^{\CC}(g)\psi(x)\label{univers2}
\end{equation}
for every $x\in X$, where $g\in {\mathcal G}$ is sufficiently close to $e$. Since the group ${\mathcal G}$ acts on $X$ locally transitively, this implies that $\psi$ maps $X$ onto an open $G_0$-invariant subset $U\subset G.p_0$. 

To show that $\psi: X\ra U$ is locally biholomorphic, assume that $\psi$ is degenerate at some point in $X$. Then property (\ref{univers2}) implies that $\psi$ is degenerate everywhere in $X$, which contradicts the openness of $U=\psi(X)$ in $G.p_0$. Thus, $\psi$ is locally biholomorphic.

{\bf Case 2.} Assume now that ${\mathfrak g}={\mathfrak{so}}_k(\CC)$ with $k\ge 5$, $k\ne 6$. Observe that the group $G_0$ contains a non-trivial compact subgroup. Hence if $G_0$ has a fixed point in $X$, by Lemma \ref{fixedpoints} and Remark \ref{fixedpointsrem} the algebra ${\mathfrak g}$ has a faithful $n$-dimensional representation. However, ${\mathfrak{so}}_k(\CC)$ for $k\ge 5$, $k\ne 6$ does not have faithful representations in dimensions less than $k-1$ (see e.g. \cite{S}, Table 1). Thus, we have shown, as in Case 1, that the $G_0$-action on $X$ is fixed point free.

Let $Y$ be any $G_0$-orbit in $X$. The group $G$ acts on the corresponding orbit $G.p_0\subset\PP^{\ell}$ almost effectively. Indeed, if $G.p_0=\{p_0\}$, the subgroup $H_0$ is discrete, and therefore $\hbox{dim}\,Y=\hbox{dim}\,G_0=k(k-1)/2$ which is impossible since $\hbox{dim}\,Y\le n$. 

Further, by Remark \ref{codimensionofh} we have $\hbox{codim}\,H\le n$. On the other hand, if $P$ is a maximal parabolic subgroup of $G$, then $\hbox{codim}\,P\ge k-2$ (see the end of Section \ref{dimtheormax}). By Proposition \ref{parabolic} we then obtain that $k=n+2$, and $H$ is a maximal parabolic subgroup of maximal dimension in $G$. It then follows that $G/H$ is biholomorphic to ${\mathcal Q}_n$ if $k\ge 7$, and to either $\PP^3$ or ${\mathcal Q}_3$ if $k=5$. Under the equivalence  $G/H\simeq{\mathcal Q}_n$ the group $G$ acts on $G/H$ as $PSO_{n+2}(\CC)$, and if $k=5$ and $G/H$ is equivalent to $\PP^3$, then under this equivalence $G$ acts on $G/H$ as $PSp_4(\CC)$ embedded in $PSL_4(\CC)$ in the standard way.

Local biholomorphicity of $\psi$ is established as in Case 1. The proof of the proposition is complete.
\end{proof}

In the following sections we will show that the map $\psi$ arising in Proposition \ref{maintheoremIprep} is in fact a biholomorphism and thus prove Theorem \ref{maintheoremI}. We will go over all real forms ${\mathfrak g}_0$ of every classical simple complex algebra ${\mathfrak g}$, determine the orbit structure for the $G_0$-action on $G/H$ in each case, and will use this information both for proving the injectivity of $\psi$ and for explicitly describing all open connected $G_0$-invariant sets $U$. Given the detailed nature of this description, it should be possible to understand basic complex-analytic properties of $X$, e.g., its Levi-convexity, Dolbeault cohomology, etc. We will consider complex algebras of the series A and C first, and those of the series B and D last. Complete lists of real forms of such algebras are well-known (see e.g. \cite{OV}, p. 233).

For future applications we restate Proposition \ref{maintheoremIprep} in more detail as follows.

\begin{proposition}\label{maintheoremIprepdetailed} Let $X$ be a connected complex manifold of dimension $n\ge 2$ admitting an almost effective action by holomorphic transformations of a connected simple Lie group $G_0$. Let ${\mathfrak g}_0$ be the Lie algebra of $G_0$, and assume that ${\mathfrak g}_0$ is of type I. Suppose, furthermore, that $k({\mathfrak g}_0)$ and $n$ satisfy Condition $(>)$. Then $k({\mathfrak g}_0)$ and $n$ satisfy Condition (=), the $G_0$-action on $X$ is fixed point free, and for every $G_0$-orbit $Y$ in $X$ the following holds:

(i) if ${\mathfrak g}={\mathfrak{sl}}_{n+1}(\CC)$, then $G/H$ is biholomorphic to $\PP^n$, and under this equivalence $G$ acts on $G/H$ as $PSL_{n+1}(\CC)$;

(ii) if ${\mathfrak g}={\mathfrak{sp}}_{n+1}(\CC)$, where $n$ is odd, then $G/H$ is biholomorphic to $\PP^n$, and under this equivalence $G$ acts on $G/H$ as $PSp_{n+1}(\CC)$ embedded in $PSL_{n+1}(\CC)$ in the standard way;

(iii) if ${\mathfrak g}={\mathfrak{so}}_{n+2}(\CC)$ for $n\ge 5$, then $G/H$ is biholomorphic to ${\mathcal Q}_n$, and under this equivalence $G$ acts on $G/H$ as $PSO_{n+2}(\CC)$; 

(iv) if ${\mathfrak g}={\mathfrak{so}}_5(\CC)$ (here $n=3$), then $G/H$ is biholomorphic to either $\PP^3$ or ${\mathcal Q}_3$; under the equivalence  $G/H\simeq\PP^3$ the group $G$ acts on $G/H$ as $PSp_4(\CC)$ embedded in $PSL_4(\CC)$ in the standard way; under the equivalence  $G/H\simeq{\mathcal Q}_3$ the group $G$ acts on $G/H$ as $SO_5(\CC)$;

(v) $S$ is open and the map $\psi_{\hbox{\tiny $S$}}$ extends to a $G_0$-equivariant locally biholomorphic surjective map $\psi: X\ra U$, where $U$ is an open $G_0$-invariant subset of $G.p_0\simeq G/H$. Here $G_0$ acts on $G/H$ as a real form of one of $PSL_{n+1}(\CC)$, $PSp_{n+1}(\CC)$ if $G/H\simeq\PP^n$, and as a real form of $PSO_{n+2}(\CC)$ if $G/H\simeq{\mathcal Q}_n$. 
\end{proposition}

\section{Classification for the A-series}\label{aseries}
\setcounter{equation}{0}

In this section we set ${\mathfrak g}={\mathfrak{sl}}_{n+1}(\CC)$. Every real form of ${\mathfrak g}$ is isomorphic to one of the algebras: ${\mathfrak{su}}_{p,\,q}$ with $p+q=n+1$, $p\ge q\ge 0$ (here ${\mathfrak{su}}_{n+1,\,0}:={\mathfrak{su}}_{n+1}$); ${\mathfrak{sl}}_{n+1}(\RR)$; ${\mathfrak{sl}}_{\hbox{\tiny$\frac{n+1}{2}$}}(\HH)$, where $\HH$ is the algebra of quaternions and in the last case $n$ is assumed to be odd. We will now separately consider each of the real forms (everywhere below we assume that $n\ge 2$).

\subsection{} Suppose first that ${\mathfrak g}_0={\mathfrak{su}}_{p,\,q}$ with $p+q=n+1$, $p\ge q\ge 0$. We start with an example. 

\begin{example}\label{exsupq} Let us realize $SU_{p,\,q}$ as the group of matrices with complex entries having determinant 1 and preserving the Hermitian form
\begin{equation}
\displaystyle\langle z,z\rangle_{p,q}:=\sum_{j=1}^{p}|z_j|^2-\sum_{j=p+1}^{n+1}|z_j|^2,\label{anglezzpq}
\end{equation}
where $z:=(z_1,\dots,z_{n+1})$ is a vector in $\CC^{n+1}$. The algebra ${\mathfrak{su}}_{p,q}$ is then realized as the Lie algebra of $SU_{p,\,q}$. The group $PSU_{p,\,q}:=SU_{p,q}/\hbox{(center)}$ acts effectively by holomorphic transformations on the projective space $\PP^n$. The orbits of this action are given by the simply-connected sets
$$
\begin{array}{lll}
\BB_{p,\,q}^{+}&:=&\left\{\zeta\in\PP^n: \langle\zeta,\zeta\rangle_{p,q}>0\right\},\\
\vspace{-0.3cm}\\
\BB_{p,\,q}^{-}&:=&\left\{\zeta\in\PP^n: \langle\zeta,\zeta\rangle_{p,q}<0\right\},\\
\vspace{-0.3cm}\\
Q_{p,\,q}&:=&\left\{\zeta\in\PP^n:\displaystyle\langle\zeta,\zeta\rangle_{p,q}=0\right\},
\end{array}
$$
where $\zeta:=(\zeta_1:\dots:\zeta_{n+1})$ are homogeneous coordinates. The sets $\BB_{p,\,q}^{+}$ and $\BB_{p,\,q}^{-}$ give open orbits. Observe that $\BB_{n+1,\,0}^{+}=\PP^n$, $\BB_{n+1,\,0}^{-}=\emptyset$, $Q_{n+1,\,0}=\emptyset$ (i.e. $PSU_{n+1}=PSU_{n+1,\,0}$ acts transitively on $\PP^n$). For $q\ge 1$ we have $\BB_{p,\,q}^{-}\ne\emptyset$, and the Levi non-degenerate hyperquadric $Q_{p,\,q}$ is a closed orbit. Observe also that $\BB_{p,p}^{+}$ is biholomorphic to $\BB_{p,p}^{-}$ for $n=2p-1$, and that $\BB_{n,\,1}^{-}=\BB^n$ (the unit ball in $\CC^n$).

We remark here for future reference that the simple connectedness of $\BB_{p,\,q}^{+}$ and $\BB_{p,\,q}^{-}$ is a manifestation of the general fact that an open orbit of a real form ${\mathfrak G}_0$ of a semisimple complex group ${\mathfrak G}$ in a complex flag manifold ${\mathfrak G}/{\mathfrak P}$ (here ${\mathfrak P}$ is a parabolic subgroup of ${\mathfrak G}$) is simply-connected (see \cite{Wo1}, Theorem 5.4 and \cite{FHW}, Proposition 4.3.5). Observe also that lower-dimensional orbits need not be simply-connected.  
\end{example}

We will now obtain our first classification result.  
  
\begin{theorem}\label{supq} Let $X$ be a connected complex manifold of dimension $n\ge 2$ admitting an almost effective action by holomorphic transformations of a connected Lie group $G_0$ with Lie algebra ${\mathfrak{su}}_{p,\,q}$, where $p+q=n+1$, $p\ge q\ge 0$. Then $X$ is biholomorphic to one of $\PP^n$, $\BB_{p,\,q}^{+}$, $\BB_{p,\,q}^{-}$.
\end{theorem}

\begin{proof} Choose $Y$ to be a $G_0$-orbit of the smallest dimension in $X$. By Proposition \ref{maintheoremIprepdetailed} the quotient $G/H$ is biholomorphic to $\PP^n$, and $G$ acts on $\PP^n$ as the group $PSL_{n+1}(\CC)$. Further, any two isomorphic real forms ${\mathfrak s}_1$, ${\mathfrak s}_2$ of a complex simple Lie algebra ${\mathfrak s}$ are obtained from one another by an inner automorphism of ${\mathfrak s}$, except in the following cases: (i) ${\mathfrak s}={\mathfrak{so}}_8(\CC)$, and ${\mathfrak s}_1$, ${\mathfrak s}_2$ are isomorphic to one of ${\mathfrak{so}}_{5,3}$, ${\mathfrak{so}}_{6,2}$, ${\mathfrak{so}}_{7,1}$; (ii) ${\mathfrak s}={\mathfrak{so}}_{4l}(\CC)$ for $l\ge 3$, and ${\mathfrak s}_1$, ${\mathfrak s}_2$ are isomorphic to ${\mathfrak{so}}_{4l}^*$ (see e.g. \cite{Dj} and a related result in \cite{Da}). It then follows that $G_0$ acts on $\PP^n$ as a subgroup of $PSL_{n+1}(\CC)$ conjugate to $PSU_{p,\,q}$. Without loss of generality, we assume that $G_0$ acts on $\PP^n$ as $PSU_{p,\,q}$. Since the map $\psi$ is locally biholomorphic, $Y$ is either an open subset or has real codimension one in $X$. We will consider these two cases separately.

Assume first that $Y$ is open in $X$, i.e. $Y=X$. In this case $U=G_0.p_0$, and therefore $U$ is either $\BB_{p,\,q}^{+}$ or $\BB_{p,\,q}^{-}$. Since $\psi:X\ra U$ is a covering map and $U$ is simply-connected, it follows that $X$ is biholomorphic to either $\BB_{p,\,q}^{+}$ or $\BB_{p,\,q}^{-}$.

Assume next that $Y$ has real codimension one in $X$. In this case $q\ge 1$, and $U$ is a $G_0$-invariant open subset of $\PP^n$ containing $Q_{p,q}$. Hence we have $U=\PP^n$. Furthermore, $\psi$ maps every codimension one $G_0$-orbit in $X$ onto $Q_{p,q}$ and every open $G_0$-orbit in $X$ onto either $\BB_{p,\,q}^{+}$ or $\BB_{p,\,q}^{-}$. The simply-connectedness of $G_0$-orbits in $\PP^n$ now yields that $\psi$ is in fact 1-to-1 on every $G_0$-orbit in $X$.

Since $\psi$ is locally biholomorphic, one can find open $G_0$-orbits $O_1$, $O_2$ such that $T:=O_1\cup O_2\cup Y$ is an open connected subset of $X$. Clearly, $\psi$ is 1-to-1 on $T$, and $\psi(T)=\PP^n$. Hence $T$ is biholomorphic to $\PP^n$, therefore $T$ is closed, which yields that $T=X$. Thus, we have shown that $X$ is biholomorphic to $\PP^n$, and the proof is complete. 
\end{proof}

Observe that Theorem \ref{supq} yields Theorem \ref{ballcharthm} stated in the introduction. The original proof of Theorem \ref{ballcharthm} in \cite{I1} was based on an explicit classification of complex $n$-dimensional manifolds admitting an effective action of the unitary group $U_n$ (see \cite{IK1}). 

\subsection{} Suppose next that ${\mathfrak g}_0={\mathfrak{sl}}_{n+1}(\RR)$. In this case all possibilities for $X$ are given by the following theorem.

\begin{theorem}\label{slr} Let $X$ be a connected complex manifold of dimension $n\ge 2$ admitting an almost effective action by holomorphic transformations of a connected Lie group $G_0$ with Lie algebra ${\mathfrak{sl}}_{n+1}(\RR)$. Then $X$ is biholomorphic to either $\PP^n$ or $\PP^n\setminus\RR\PP^n$, where $\RR\PP^n$ is the real projective space in $\PP^n$.
\end{theorem}

\begin{proof} As before, let $Y$ be a $G_0$-orbit of the smallest dimension in $X$. As in the proof of Theorem \ref{supq}, we see that $G_0$ acts on $G/H\simeq\PP^n$ as a subgroup of $PSL_{n+1}(\CC)$ conjugate to $PSL_{n+1}(\RR)$, and without loss of generality we assume that $G_0$ acts as $PSL_{n+1}(\RR)$. Observe that for $n\ge 2$ the group $PSL_{n+1}(\RR)$ has exactly two orbits in $\PP^n$: the totally real closed orbit $\RR\PP^n$ and the open orbit $\PP^n\setminus\RR\PP^n$ (see \cite{Wi}, p. 209). Since $\psi$ is locally biholomorphic, $Y$ is either an open subset or a totally real $n$-dimensional submanifold in $X$.  

Assume first that $Y$ is open in $X$. In this case $U=G_0.p_0$, and therefore $U$ coincides with $\PP^n\setminus\RR\PP^n$. Using the simple connectedness of $U$ we now see that $X$ is biholomorphic to $\PP^n\setminus\RR\PP^n$.

Assume next that $Y$ is totally real in $X$. In this case $U$ is a $G_0$-invariant open subset of $\PP^n$ containing $\RR\PP^n$. Hence we have $U=\PP^n$. Furthermore, $\psi$ maps every totally real $G_0$-orbit in $X$ onto $\RR\PP^n$ and every open $G_0$-orbit in $X$ onto $\PP^n\setminus\RR\PP^n$. The map $\psi$ is 1-to-1 on every open $G_0$-orbit in $X$ and is either a 2-to-1 covering map or a diffeomorphism on every totally real orbit. We will now show that $\psi$ is in fact 1-to-1 on every totally real orbit. Let $Y'$ be such an orbit and assume that for two distinct points $x_1,x_2\in Y'$ we have $\psi(x_1)=\psi(x_2)$. Choose non-intersecting neighborhoods $V_j$ of $x_j$, such that $\psi$ is 1-to-1 on $V_j$, $j=1,2$, and $\psi(V_1)=\psi(V_2):=W$. It then follows that there exists an open $G_0$-orbit $O'$ in $X$ such that $W_j:=V_j\setminus(V_j\cap Y')$ is a connected open set lying in $O'$ for each $j$. Clearly, we have $\psi(W_1)=\psi(W_2)=W\setminus(W\cap\RR\PP^n)$, which contradicts the fact that $\psi$ is 1-to-1 on $O'$. Thus, we have shown that $\psi$ is 1-to-1 on every $G_0$-orbit in $X$.

Next, one can find an open $G_0$-orbit $O$ such that $T:=O\cup Y$ is an open connected subset of $X$. Clearly, $\psi$ is 1-to-1 on $T$, and $\psi(T)=\PP^n$. Hence $T$ is biholomorphic to $\PP^n$, therefore $T$ is closed, which implies that $T=X$. Thus, $X$ is biholomorphic to $\PP^n$, and the proof is complete. 
\end{proof}

\subsection{} Suppose finally that $n\ge 3$ is odd and ${\mathfrak g}_0={\mathfrak{sl}}_{\hbox{\tiny$\frac{n+1}{2}$}}(\HH)$. It turns out that in this case there is only one possibility for $X$.  

\begin{theorem}\label{slh} Let $X$ be a connected complex manifold of odd dimension $n\ge 3$ admitting an almost effective action by holomorphic transformations of a connected Lie group $G_0$ with Lie algebra ${\mathfrak{sl}}_{{\hbox{\tiny$\frac{n+1}{2}$}}}(\HH)$. Then $X$ is biholomorphic to $\PP^n$.
\end{theorem}

\begin{proof} Choose $Y$ to be any $G_0$-orbit in $X$. It follows as above that $G_0$ acts on $G/H\simeq\PP^n$ as a subgroup of $PSL_{n+1}(\CC)$ conjugate to $PSL_{{\hbox{\tiny$\frac{n+1}{2}$}}}(\HH)$. Observe that the group $PSL_{{\hbox{\tiny$\frac{n+1}{2}$}}}(\HH)$ acts transitively on $\PP^n$ (see \cite{Wi}, p. 211 and \cite{Wo2}, Corollary 1.7), and therefore $Y$ is open in $X$. Thus, $Y=X$ and $\psi:X\ra\PP^n$ is a covering map. Hence $X$ is biholomorphic to $\PP^n$.
\end{proof}

Theorems \ref{supq}, \ref{slr}, \ref{slh} complete the proof of Theorem \ref{maintheoremI} in the case when ${\mathfrak g}$ belongs to the A-series.

\section{Classification for the C-series}\label{cseries}
\setcounter{equation}{0}

In this section we assume that $n$ is odd and set ${\mathfrak g}={\mathfrak{sp}}_{n+1}(\CC)$ with $n+1=2l$, $l\ge 3$. Every real form of ${\mathfrak g}$ is isomorphic to one of the algebras: ${\mathfrak{sp}}_{p,\,q}$ with $p+q=l$, $p\ge q\ge 0$ (here ${\mathfrak{sp}}_{l,\,0}:={\mathfrak{sp}}_l$) ; ${\mathfrak{sp}}_{n+1}(\RR)$. We will now separately consider each of the real forms (everywhere below we assume that $n\ge 5$ is odd).

\subsection{} Suppose first that ${\mathfrak g}_0={\mathfrak{sp}}_{p,\,q}$ with $p+q=(n+1)/2$, $p\ge q\ge 0$. We start with an example.

\begin{example}\label{exsppq} As before, let us realize $Sp_{n+1}(\CC)$ as the group of matrices with complex entries preserving the skew-symmetric form given by (\ref{matrixr}) for $l=(n+1)/2$. Next, realize the group $Sp_{p,q}$ as the real subgroup of $Sp_{n+1}(\CC)$ that consists of matrices preserving the Hermitian form
$$
(z,z)_{p,q}:=\sum_{j=1}^{p}|z_j|^2-\sum_{j=p+1}^l|z_j|^2+\sum_{j=l+1}^{l+p}|z_j|^2-\sum_{j=l+p+1}^{n+1}|z_j|^2,
$$
where $z:=(z_1,\dots,z_{n+1})$ is a vector in $\CC^{n+1}$. The algebra ${\mathfrak{sp}}_{p,q}$ is then realized as the Lie algebra of $Sp_{p,q}$. The group $PSp_{p,\,q}:=Sp_{p,\,q}/\hbox{(center)}$ acts by holomorphic transformations on the projective space $\PP^n$. Its orbits are given by the simply-connected sets
$$
\begin{array}{lll}
\hat\BB_{p,\,q}^{+}&:=&\left\{\zeta\in\PP^n: (\zeta,\zeta)_{p,q}>0\right\},\\
\vspace{-0.3cm}\\
\hat\BB_{p,\,q}^{-}&:=&\left\{\zeta\in\PP^n: (\zeta,\zeta)_{p,q}<0\right\},\\
\vspace{-0.3cm}\\
\hat Q_{p,\,q}&:=&\left\{\zeta\in\PP^n:\displaystyle (\zeta,\zeta)_{p,q}=0\right\},
\end{array}
$$
where $\zeta:=(\zeta_1:\dots:\zeta_{n+1})$ are homogeneous coordinates. Clearly, there exists an automorphism of $\PP^n$ that maps $\hat\BB_{p,\,q}^{+}$, $\hat\BB_{p,\,q}^{-}$, $\hat Q_{p,\,q}$ onto the sets $\BB_{2p,\,2q}^{+}$, $\BB_{2p,\,2q}^{-}$, $Q_{2p,\,2q}$, respectively (see Example \ref{exsupq}). The sets $\hat\BB_{p,\,q}^{+}$ and $\hat\BB_{p,\,q}^{-}$ give open orbits. For $q\ge 1$ we have $\hat\BB_{p,\,q}^{-}\ne\emptyset$, and the Levi non-degenerate hyperquadric $\hat Q_{p,\,q}$ is the closed orbit (see \cite{Wi}, p. 221). 
\end{example}

The proof of the following theorem is completely analogous to that of Theorem \ref{supq}, and we omit it.

\begin{theorem}\label{sppq} Let $X$ be a connected complex manifold of odd dimension $n\ge 5$ admitting an almost effective action by holomorphic transformations of a connected Lie group $G_0$ with Lie algebra ${\mathfrak{sp}}_{p,\,q}$, where $p+q=(n+1)/2$, $p\ge q\ge 0$. Then $X$ is biholomorphic to one of $\PP^n$, $\BB_{2p,\,2q}^{+}$, $\BB_{2p,\,2q}^{-}$.
\end{theorem}

\subsection{} Suppose next that ${\mathfrak g}_0={\mathfrak{sp}}_{n+1}(\RR)$. We start with an example.

\begin{example}\label{exspr} We realize the group $Sp_{n+1}(\RR)$ as $Sp_{n+1}(\CC)\cap GL_{2(n+1)}(\RR)$ and the algebra ${\mathfrak{sp}}_{n+1}(\RR)$ as the Lie algebra of $Sp_{n+1}(\RR)$. The group $PSp_{n+1}(\RR):=Sp_{n+1}(\RR)/\hbox{(center)}$ acts effectively on $\PP^n$, and $\RR\PP^n\subset\PP^n$ is clearly a closed $PSp_{n+1}(\RR)$-orbit. We will now describe the remaining orbits. For a pair of vectors $z,w\in\CC^{n+1}$ consider the skew-symmetric form $(z,w)$ given by (\ref{matrixr}) for $l=(n+1)/2$. Define
$$
\begin{array}{lll}
D_n^{+}&:=&\left\{\zeta\in\PP^n: (\hbox{Re}\,\zeta,\hbox{Im}\,\zeta)>0\right\},\\
\vspace{-0.3cm}\\
D_n^{-}&:=&\left\{\zeta\in\PP^n: (\hbox{Re}\,\zeta,\hbox{Im}\,\zeta)<0\right\},\\
\vspace{-0.3cm}\\
\Sigma_n&:=&\left\{\zeta\in\PP^n: (\hbox{Re}\,\zeta,\hbox{Im}\,\zeta)=0,\,\,\hbox{Re}\,\zeta\ne0,\,\,  \hbox{Im}\,\zeta\ne0 \right\},
\end{array}
$$
where $\zeta:=(\zeta_1:\dots:\zeta_{n+1})$ are homogeneous coordinates. The sets $D_n^{+}$ and $D_n^{-}$ are open simply-connected $PSp_{n+1}(\RR)$-orbits, and $\Sigma_n$ is a $PSp_{n+1}(\RR)$-orbit of real codimension one (see \cite{Wi}, p. 219). Note that there exists a holomorphic automorphism of $\PP^n$ that maps $\Sigma_n\cup\RR\PP^n$ onto the hyperquadric $Q_{{\hbox{\tiny$\frac{n+1}{2}$}},\,{\hbox{\tiny$\frac{n+1}{2}$}}}$ introduced in Example \ref{exsupq} (see \cite{Wi}, p. 220). In particular, each of $D_n^{+}$, $D_n^{-}$ is biholomorphic to $\BB_{{\hbox{\tiny$\frac{n+1}{2}$}},\,{\hbox{\tiny$\frac{n+1}{2}$}}}^{+}$.   
\end{example}

The result of this subsection is the following theorem. 

\begin{theorem}\label{spr} Let $X$ be a connected complex manifold of odd dimension $n\ge 5$ admitting an almost effective action by holomorphic transformations of a connected Lie group $G_0$ with Lie algebra ${\mathfrak{sp}}_{n+1}(\RR)$. Then $X$ is biholomorphic to one of $\PP^n$, $\PP^n\setminus\RR\PP^n$, $\BB_{{\hbox{\tiny$\frac{n+1}{2}$}},\,{\hbox{\tiny$\frac{n+1}{2}$}}}^{+}$.
\end{theorem}

\begin{proof} Choose $Y$ to be a $G_0$-orbit of the smallest dimension in $X$. By Proposition \ref{maintheoremIprepdetailed} the quotient $G/H$ is biholomorphic to $\PP^n$, and $G$ acts on $\PP^n$ as the group $PSp_{n+1}(\CC)$ embedded in $PSL_{n+1}(\CC)$ in the standard way. As in the proof of Theorem \ref{supq}, we see that $G_0$ acts on $\PP^n$ as a subgroup of $PSp_{n+1}(\CC)$ conjugate to $PSp_{n+1}(\RR)$. Without loss of generality, we assume that $G_0$ acts on $\PP^n$ as $PSp_{n+1}(\RR)$. Since the map $\psi$ is locally biholomorphic, $Y$ is either (i) an open subset, or (ii) a real hypersurface, or (iii) a totally real $n$-dimensional submanifold. We will consider these three cases separately.

The case when $Y$ is open is treated as the corresponding case in the proof of Theorem \ref{supq}, with the result that $X$ is biholomorphic to one of $D_n^{+}$, $D_n^{-}$ and hence to $\BB_{{\hbox{\tiny$\frac{n+1}{2}$}},\,{\hbox{\tiny$\frac{n+1}{2}$}}}^{+}$. 

Assume next that $Y$ has real codimension one in $X$. In this case $U$ is a $G_0$-invariant open subset in $\PP^n$ containing $\Sigma_n$ and not containing $\RR\PP^n$. Hence we have $U=\PP^n\setminus\RR\PP^n$. Clearly, $\psi$ maps every open $G_0$-orbit in $X$ biholomorphically onto one of the two open $G_0$-orbits $D_n^{+}$, $D_n^{-}$. Furthermore, every $G_0$-orbit of codimension one covers $\Sigma_n$ by means of $\psi$. We will now show that $\psi$ is in fact 1-to-1 on every codimension one orbit. Let $Y'$ be such an orbit and assume that for two distinct points $x_1,x_2\in Y'$ we have $\psi(x_1)=\psi(x_2)$. Choose non-intersecting neighborhoods $V_j$ of $x_j$, such that $\psi$ is 1-to-1 on $V_j$, $j=1,2$, and $\psi(V_1)=\psi(V_2):=W$. It then follows that there exist open $G_0$-orbits $O'$ and $O''$ in $X$ such that $W_j:=V_j\setminus(V_j\cap Y')$ is an open set having exactly two connected components $W_j'$, $W_j''$, with $W_j'$ lying in $O'$ and $W_j''$ lying in $O''$ for each $j$. Clearly, we have either $\psi(W_1')=\psi(W_2')=W\cap D_n^{+}$ and $\psi(W_1'')=\psi(W_2'')=W\cap D_n^{-}$, or $\psi(W_1')=\psi(W_2')=W\cap D_n^{-}$ and $\psi(W_1'')=\psi(W_2'')=W\cap D_n^{+}$. However, this contradicts the fact that $\psi$ is 1-to-1 on each of $O'$, $O''$. Thus, we have shown that $\psi$ is 1-to-1 on every $G_0$-orbit in $X$.

As in the proof of Theorem \ref{supq}, one can find open $G_0$-orbits $O_1$, $O_2$ such that $T:=O_1\cup O_2\cup Y$ is an open connected subset of $X$. Clearly, $\psi$ is 1-to-1 on $T$ and maps it onto $\PP^n\setminus\RR\PP^n$. Suppose that $T\ne X$ and consider a point $x\in\partial T$. Let $Y''$ be the $G_0$-orbit of $x$. Clearly, $Y''$ has codimension one. Further, since $\psi$ is locally biholomorphic and $\PP^n$ near every point of $\Sigma_n$ splits into the disjoint union of a portion of $\Sigma_n$ and two non-empty open subsets lying in the distinct open orbits $D_n^{+}$ and $D_n^{-}$, there exists a neighborhood ${\mathcal V}$ of $x$ on which $\psi$ is 1-to-1 and such that ${\mathcal V}={\mathcal V}_1\cup {\mathcal V}_2\cup ({\mathcal V}\cap Y'')$, where ${\mathcal V}_1$ is the intersection of ${\mathcal V}$ with an open orbit $O'''$ contained in $T$ (i.e. $O'''$ is one of $O_1$, $O_2$), and ${\mathcal V}_2$ is the intersection of ${\mathcal V}$ with some other open $G_0$-orbit in $X$. Since each of $Y$ and $Y''$ is mapped by $\psi$ onto $\Sigma_n$, there exists a point $z\in Y$ with $\psi(z)=\psi(x)$. Let ${\mathcal U}$ be a neighborhood of $z$ not intersecting ${\mathcal V}$, on which $\psi$ is 1-to-1 and such that $\psi({\mathcal U})\subset\psi({\mathcal V})$. Then $\psi({\mathcal U}\cap O''')\subset\psi({\mathcal V}_1)$, which implies that $\psi$ is not 1-to-1 on $O'''$. This contradiction shows that in fact we have $T=X$, and thus $X$ is biholomorphic to $\PP^n\setminus\RR\PP^n$.

Assume finally that $Y$ is totally real in $X$. In this case $U$ is a $G_0$-invariant open subset of $\PP^n$ containing $\RR\PP^n$. Hence we have $U=\PP^n$. In particular, there exists a codimension one orbit in $X$. Clearly, $\psi$ maps every totally real $G_0$-orbit in $X$ onto $\RR\PP^n$. Since $\psi$ is locally biholomorphic, the union of totally real orbits in $X$ is closed and its complement $X'$ is connected. Arguing as above, we obtain that $\psi$ maps $X'$ biholomorphically onto $\PP^n\setminus\RR\PP^n$. Fix $y\in\partial X'$ and consider its $G_0$-orbit $Y'''$ (clearly, $Y'''$ is totally real). We will now show that $\psi$ is 1-to-1 on $Y'''$. Assume that for two distinct points $y_1,y_2\in Y'''$ we have $\psi(y_1)=\psi(y_2)$. Choose non-intersecting neighborhoods ${\mathfrak V}_j$ of $y_j$, such that $\psi$ is 1-to-1 on ${\mathfrak V}_j$, $j=1,2$, and $\psi({\mathfrak V}_1)=\psi({\mathfrak V}_2):={\mathfrak W}$. It then follows that there exist open $G_0$-orbits $\hat O$ and $\tilde O$ and a codimension one orbit $Y_0$ in $X$ such that ${\mathfrak W}_j:={\mathfrak V}_j\setminus({\mathfrak V}_j\cap (Y'''\cup Y_0))$ is an open set having exactly two connected components $\hat {\mathfrak W}_j$, $\tilde{\mathfrak W}_j$, with $\hat {\mathfrak W}_j$ lying in $\hat O$ and $\tilde {\mathfrak W}_j$ lying in $\tilde O$ for each $j$. Clearly, we have either $\psi(\hat{\mathfrak W}_1)=\psi(\hat{\mathfrak W}_2)={\mathfrak W}\cap D_n^{+}$ and $\psi(\tilde{\mathfrak W}_1)=\psi(\tilde{\mathfrak W}_2)={\mathfrak W}\cap D_n^{-}$, or $\psi(\hat{\mathfrak W}_1)=\psi(\hat{\mathfrak W}_2)={\mathfrak W}\cap D_n^{-}$ and $\psi(\tilde{\mathfrak W}_1)=\psi(\tilde{\mathfrak W}_2)={\mathfrak W}\cap D_n^{+}$. However, this contradicts the fact that $\psi$ is 1-to-1 on each of $\hat O$, $\tilde O$. Thus, we have shown that $\psi$ is 1-to-1 on $Y'''$ (and in fact on every totally real $G_0$-orbit in $X$).

Clearly, $T:=X'\cup Y'''$ is a connected open subset of $X$ that $\psi$ maps biholomorphically onto $\PP^n$. Therefore, $T$ is also closed, which implies that $T=X$. Thus, we have shown that $X$ is biholomorphic to $\PP^n$. The proof is complete.\end{proof}

Theorems \ref{sppq}, \ref{spr} complete the proof of Theorem \ref{maintheoremI} in the case when ${\mathfrak g}$ belongs to the C-series.

\section{Classification for the B- and D-series}\label{bseries}
\setcounter{equation}{0}

In this section we set ${\mathfrak g}={\mathfrak{so}}_{n+2}(\CC)$ and assume that $n\ge 3$, $n\ne 4$. Every real form of ${\mathfrak g}$ is isomorphic to one of the algebras: ${\mathfrak{so}}_{p,\,q}$ with $p+q=n+2$, $p\ge q\ge 0$ (here ${\mathfrak{so}}_{n+2,\,0}:={\mathfrak{so}}_{n+2}(\RR)$); ${\mathfrak{so}}^*_{n+2}$, where in the last case $n\ge 8$ is even (recall that ${\mathfrak{so}}^*_8$ is isomorphic to ${\mathfrak{so}}_{6,2}$). We will now separately consider each of the real forms.

\subsection{} Suppose first that ${\mathfrak g}_0={\mathfrak{so}}_{p,\,q}$ with $p+q=n+2$, $p\ge q\ge 0$, and that for $n=6$ the algebra ${\mathfrak g}_0$ is not one of ${\mathfrak{so}}_{5,\,3}$, ${\mathfrak{so}}_{6,\,2}$, ${\mathfrak{so}}_{7,\,1}$. We start with an example. 

\begin{example}\label{exsopq} Fix $p, q$ and consider the group ${\mathcal{SO}}_{p,q}$ of $(n+2)\times(n+2)$-matrices with real entries having determinant 1 and preserving the quadratic form
\begin{equation}
\displaystyle [z,z]_{p,q}:=\sum_{j=1}^{p}z_j^2-\sum_{j=p+1}^{n+2}z_j^2,\label{zzpq}
\end{equation}
where $z:=(z_1,\dots,z_{n+2})$ is a vector in $\CC^{n+2}$. Let
$$
{\mathcal A}_{p,q}:=\left(
\begin{array}{cc}
I_p & 0\\
0 & iI_q
\end{array}
\right),
$$
where for any $m$ we denote by $I_m$ the $m\times m$ identity matrix. We now realize the group $SO_{p,\,q}$ as the subgroup ${\mathcal A}_{p,q}{\mathcal{SO}}_{p,q}{\mathcal A}_{p,q}^{-1}$ of $SO_{n+2}(\CC)$ and the algebra ${\mathfrak{so}}_{p,q}$ as the Lie algebra of $SO_{p,\,q}$. The group $PSO_{p,\,q}^{\circ}:=SO_{p,\,q}^{\circ}/\hbox{(center)}$ acts by holomorphic transformations on the projective quadric ${\mathcal Q}_n$. The orbits of this action are described by means of the following sets
$$
\begin{array}{lll}
\Omega_{p,q}^{+}&:=&\left\{\zeta\in{\mathcal Q}_n: \left[\hbox{Re}\,\left({\mathcal A}_{p,q}^{-1}\zeta\right),\hbox{Re}\,\left({\mathcal A}_{p,q}^{-1}\zeta\right)\right]_{p,q}>0\right\},\\
\vspace{-0.3cm}\\
\Omega_{p,q}^{-}&:=&\left\{\zeta\in{\mathcal Q}_n: \left[\hbox{Re}\,\left({\mathcal A}_{p,q}^{-1}\zeta\right),\hbox{Re}\,\left({\mathcal A}_{p,q}^{-1}\zeta\right)\right]_{p,q}<0\right\},\\
\vspace{-0.3cm}\\
{\mathcal S}_{p,q}&:=&\left\{\zeta\in{\mathcal Q}_n: \left[\hbox{Re}\,\left({\mathcal A}_{p,q}^{-1}\zeta\right),\hbox{Re}\,\left({\mathcal A}_{p,q}^{-1}\zeta\right)\right]_{p,q}=0\right\}.
\end{array}
$$
The sets $\Omega_{p,q}^{+}$ and $\Omega_{p,q}^{-}$ give open orbits. They are are simply-connected (see \cite{Wo1}, Theorem 5.4 and \cite{FHW}, Proposition 4.3.5). Observe that $\Omega_{n+2,0}^{+}={\mathcal Q}_n$, $\Omega_{n+2,0}^{-}=\emptyset$, ${\mathcal S}_{n+2,0}=\emptyset$ (i.e. $PSO_{n+2}^{\circ}(\RR)=PSO_{n+2,0}^{\circ}$ acts transitively on ${\mathcal Q}_n$). Further, $\Omega_{n+1,1}^{-}=\emptyset$, and ${\mathcal S}_{n+1,1}$ is a closed totally real $n$-dimensional $PSO_{n+1,1}^{\circ}$-orbit. For $q\ge 2$ we have $\Omega_{p,q}^{-}\ne\emptyset$, and ${\mathcal S}_{p,q}={\mathcal S}_{p,q}^1\cup {\mathcal S}_{p,q}^2$, where ${\mathcal S}_{p,q}^1$ is a closed totally real $n$-dimensional orbit and ${\mathcal S}_{p,q}^2$ is an orbit of real codimension one (${\mathcal S}_{p,q}^1$ and ${\mathcal S}_{p,q}^2$ are the orbits of the points ${\mathcal A}_{p,q}\bigl((1:0:\dots:0:1)\bigr)$ and ${\mathcal A}_{p,q}\bigl((1:i:0:\dots:0:i:1)\bigr)$, respectively). For convenience we set ${\mathcal S}_{n+1,1}^1:={\mathcal S}_{n+1,1}$ and ${\mathcal S}_{n+1,1}^2:=\emptyset$. Observe that the closed orbit ${\mathcal S}_{p,q}^1$ is diffeomorphic to the real projective quadric
$$
\left\{x\in\RR\PP^{n+1}: [x,x]_{p,q}=0\right\},
$$
where $x:=(x_1:\dots:x_{n+2})$ are homogeneous coordinates in $\RR\PP^{n+1}$.
\end{example}

The first result of this section is the following theorem.

\begin{theorem}\label{sopq} Let $X$ be a connected complex manifold of dimension $n\ge 5$ admitting an almost effective action by holomorphic transformations of a connected Lie group $G_0$ with Lie algebra ${\mathfrak{so}}_{p,\,q}$, where $p+q=n+2$, $p\ge q\ge 0$. Assume that for $n=6$ we have $p\ne 5,6,7$. Then $X$ is biholomorphic to one of ${\mathcal Q}_n$, ${\mathcal Q}_n\setminus{\mathcal S}_{p,q}^1$, $\Omega_{p,\,q}^{+}$, $\Omega_{p,\,q}^{-}$.
\end{theorem}

\begin{proof} Choose $Y$ to be a $G_0$-orbit of the smallest dimension in $X$. By Proposition \ref{maintheoremIprepdetailed} the quotient $G/H$ is biholomorphic to ${\mathcal Q}_n$, and $G$ acts on ${\mathcal Q}_n$ as the group $PSO_{n+2}(\CC)$. As in the proof of Theorem \ref{supq}, we see that $G_0$ acts on ${\mathcal Q}_n$ as a subgroup of $PSO_{n+2}(\CC)$ conjugate to $PSO_{p,q}^{\circ}$. Since the map $\psi$ is locally biholomorphic, $Y$ is either (i) an open subset, or (ii) a real hypersurface, or (iii) a totally real $n$-dimensional submanifold.

The case when $Y$ is open is treated as the corresponding case in the proof of Theorem \ref{supq}, and one obtains that $X$ is biholomorphic to one of $\Omega_{p,\,q}^{+}$, $\Omega_{p,\,q}^{-}$. 

If $Y$ is a real hypersurface in $X$ (here $q\ge 2$), we argue as in the corresponding case in the proof of Theorem \ref{spr}.

If $Y$ is totally real, then $q\ge 1$ and two situations are possible. If $q=1$, then $X$ contains no codimension one orbits, and the proof follows as for Theorem \ref{slr}. If $q\ge 2$, then $X$ contains a codimension one orbit, and we argue as in the corresponding case in the proof of Theorem \ref{spr}.
\end{proof}

We will now deal with the case $n=3$ not covered by Theorem \ref{sopq}. Note that ${\mathfrak{so}}_5(\CC)$ is isomorphic to ${\mathfrak{sp}}_4(\CC)$, and the real forms of these algebras are ${\mathfrak{so}}_{5,\,0}:={\mathfrak{so}}_5(\RR)\simeq {\mathfrak{sp}}_{2,\,0}:={\mathfrak{sp}}_2$; ${\mathfrak{so}}_{4,\,1}\simeq{\mathfrak{sp}}_{1,\,1}$; ${\mathfrak{so}}_{3,\,2}\simeq{\mathfrak{sp}}_4(\RR)$. The group $SO_5(\CC)\simeq PSp_4(\CC)$ acts on each of $\PP^3$ and ${\mathcal Q}_3$ (see \cite{Wi}, pp. 214--224 for a thorough description of these actions), and our classification result for $n=3$ is a combination of statements of Theorems \ref{sppq}, \ref{spr}, \ref{sopq}.

\begin{theorem}\label{sopqn=3} Let $X$ be a connected complex manifold of dimension 3 admitting an almost effective action by holomorphic transformations of a connected Lie group $G_0$ with Lie algebra ${\mathfrak{so}}_{p,\,q}$, where $p+q=5$, $p\ge q\ge 0$. Then $X$ is biholomorphic to one of 

(i) $\PP^3$, ${\mathcal Q}_3$, if $p=5$;

(ii) $\PP^3$, ${\mathcal Q}_3$, $\BB_{2,2}^{+}$, $\Omega_{4,1}^{+}$ , if $p=4$;

(iii) $\PP^3$, ${\mathcal Q}_3$, $\PP^3\setminus\RR\PP^3$, ${\mathcal Q}_3\setminus{\mathcal S}_{3,2}^1$, $\BB_{2,2}^{+}$, $\Omega_{3,2}^{+}$, $\Omega_{3,2}^{-}$, if $p=3$.
\end{theorem}

The proof of Theorem \ref{sopqn=3} is obtained by suitably combining the proofs of Theorems \ref{sppq}, \ref{spr}, \ref{sopq}, and we omit it. 

\subsection{} Suppose now that $n\ge 8$ is even and ${\mathfrak g}_0={\mathfrak{so}}^*_{n+2}$. We start with an example.

\begin{example}\label{exsostar} Consider the real subgroup $SO_{n+2}^*:=SO_{n+2}(\CC)\cap SL_{\hbox{\tiny$\frac{n+2}{2}$}}(\HH)$ of $SO_{n+2}(\CC)$ and realize ${\mathfrak{so}}^*_{n+2}$ as the Lie algebra of $SO_{n+2}^*$. Matrices from $SO_{n+2}^*$ preserve the skew-symmetric Hermitian form
$$
[z,z]:=\sum_{j=1}^{\frac{n+2}{2}}\left(z_j\overline{z}_{\frac{n+2}{2}+j}-z_{\frac{n+2}{2}+j}\overline{z}_j\right).
$$   
The group $PSO_{n+2}^*:=SO_{n+2}^*/\hbox{(center)}$ acts effectively on ${\mathcal Q}_n$. The orbits of this action are described as follows. Define
$$
\begin{array}{lll}
E_n^{+}&:=&\left\{\zeta\in{\mathcal Q}_n: i[\zeta,\zeta]>0\right\},\\
\vspace{-0.3cm}\\
E_n^{-}&:=&\left\{\zeta\in{\mathcal Q}_n: i[\zeta,\zeta]<0\right\},\\
\vspace{-0.3cm}\\
{\mathfrak S}_n&:=&\left\{\zeta\in{\mathcal Q}_n: [\zeta,\zeta]=0\right\},
\end{array}
$$
where $\zeta:=(\zeta_1:\dots:\zeta_{n+2})$ are homogeneous coordinates in $\PP^{n+1}$. The sets $E_n^{+}$ and $E_n^{-}$ are open simply-connected $PSO_{n+2}^*$-orbits, and ${\mathfrak S}_n$ is a closed $PSO_{n+2}^*$-orbit of real codimension one.

Next, for any $m$ consider the $m\times m$-matrix
\begin{equation}
{\mathcal R}_m:=\left(
\begin{array}{rc}
-1 & 0\\
0 & I_{m-1}
\end{array}
\right).\label{matricmathcalr}
\end{equation}
The matrix ${\mathcal R}_{n+2}$ has determinant $-1$ and lies in the connected component of $O_{n+2}(\CC)$ complementary to $SO_{n+2}(\CC)$. Conjugation by ${\mathcal R}_{n+2}$ stabilizes $SO_{n+2}(\CC)$ and induces an outer automorphism of this group. Define $\widehat{SO}_{n+2}^*$ to be the image of $SO_{n+2}^*$ under this automorphism, and let $\widehat{{\mathfrak{so}}}_{n+2}^*$ be the Lie algebra of $\widehat{SO}_{n+2}^*$. The subgroups $\widehat{SO}_{n+2}^*$ and $SO_{n+2}^*$ are conjugate to each other in $SO_{n+2}(\CC)$ only if $(n+2)/2$ is odd (see e.g. \cite{Dj}). The quotient $P\widehat{SO}_{n+2}^*:=\widehat{SO}_{n+2}^*/\hbox{(center)}$ acts effectively on ${\mathcal Q}_n$, and the automorphism of ${\mathcal Q}_n$ given by the matrix ${\mathcal R}_{n+2}$ transforms the orbits of $P\widehat{SO}_{n+2}^*$ into those of $PSO_{n+2}^*$ and vice versa.    
\end{example}

We are now ready to prove the following theorem.

\begin{theorem}\label{sostar} Let $X$ be a connected complex manifold of even dimension $n\ge 8$ admitting an almost effective action by holomorphic transformations of a connected Lie group $G_0$ with Lie algebra ${\mathfrak{so}}^*_{n+2}$. Then $X$ is biholomorphic to one of ${\mathcal Q}_n$, $E_n^{+}$, $E_n^{-}$.
\end{theorem}

\begin{proof} Choose $Y$ to be a $G_0$-orbit of the smallest dimension in $X$. By Proposition \ref{maintheoremIprepdetailed} the quotient $G/H$ is biholomorphic to ${\mathcal Q}_n$, and $G$ acts on ${\mathcal Q}_n$ as the group $PSO_{n+2}(\CC)$. Then the group $G_0$ acts on ${\mathcal Q}_n$ as a subgroup of $PSO_{n+2}(\CC)$ conjugate either to $PSO^*_{n+2}$ or, if $(n+2)/2$ is even, to $P\widehat{SO}_{n+2}^*$. Since the map $\psi$ is locally biholomorphic, $Y$ is either open or has codimension one in $X$.

We now proceed as in the proof of Theorem \ref{supq} and obtain that $X$ is biholomorphic either to one of ${\mathcal Q}_n$, $E_n^{+}$, $E_n^{-}$, or, if $(n+2)/2$ is even, to one of the open $P\widehat{SO}_{n+2}^*$-orbits. Since each of the latter is equivalent to one of  $E_n^{+}$, $E_n^{-}$, the theorem follows.  
\end{proof}

\subsection{} Suppose finally that $n=6$ and ${\mathfrak g}_0$ is one of the algebras ${\mathfrak{so}}_{5,3}$, ${\mathfrak{so}}_{6,2}$, ${\mathfrak{so}}_{7,1}$. In this case there are exactly three ways (up to inner automorphisms of ${\mathfrak{so}}_8(\CC)$) to embed ${\mathfrak g}_0$ into ${\mathfrak{so}}_8(\CC)$ as a real form (see e.g. \cite{Dj}). The equivalence classes of embeddings -- with respect to the group $\hbox{Inn}\,({\mathfrak{so}}_8(\CC))$ of inner automorphisms of ${\mathfrak{so}}_8(\CC)$ -- are described as follows. One equivalence class is represented by the standard embeddings ${\mathfrak{so}}_{5,3}^0$, ${\mathfrak{so}}_{6,2}^0$, ${\mathfrak{so}}_{7,1}^0$ of ${\mathfrak{so}}_{5,3}$, ${\mathfrak{so}}_{6,2}$, ${\mathfrak{so}}_{7,1}$ into ${\mathfrak{so}}_8(\CC)$, respectively. These are the subalgebras of ${\mathfrak{so}}_8(\CC)$ corresponding to the subgroups $SO_{5,3}^{\circ}$, $SO_{6,2}^{\circ}$, $SO_{7,1}^{\circ}$ of $SO_8(\CC)$, respectively, described in Example \ref{exsopq}. The other two equivalence classes are represented by ${\mathfrak{so}}_{5,3}^1:=\theta({\mathfrak{so}}_{5,3}^0)$, ${\mathfrak{so}}_{6,2}^1:=\theta({\mathfrak{so}}_{6,2}^0)$, ${\mathfrak{so}}_{7,1}^1:=\theta({\mathfrak{so}}_{7,1}^0)$ and by ${\mathfrak{so}}_{5,3}^2:=\theta^2({\mathfrak{so}}_{5,3}^0)$, ${\mathfrak{so}}_{6,2}^2:=\theta^2({\mathfrak{so}}_{6,2}^0)$, ${\mathfrak{so}}_{7,1}^2:=\theta^2({\mathfrak{so}}_{7,1}^0)$, where $\theta$ is a certain outer automorphism of ${\mathfrak{so}}_8(\CC)$ (a triality automorphism). The automorphism $\theta$ satisfies the condition $\theta^3\in\hbox{Inn}\,({\mathfrak{so}}_8(\CC))$ and has the form $\theta=\theta_0\circ{\varphi}$. Here $\theta_0$ is an outer automorphism of ${\mathfrak{so}}_8(\CC)$ leaving invariant the equivalence classes of ${\mathfrak{so}}_{5,3}^2$, ${\mathfrak{so}}_{6,2}^2$, ${\mathfrak{so}}_{7,1}^2$ and such that $\theta_0^2\in\hbox{Inn}\,({\mathfrak{so}}_8(\CC))$, and $\varphi$ is the outer automorphism of ${\mathfrak{so}}_8(\CC)$ given by the conjugation by the matrix ${\mathcal R}_8\in O_8(\CC)\setminus SO_8(\CC)$ (see (\ref{matricmathcalr})). Therefore, the equivalence classes of ${\mathfrak{so}}_{5,3}^2$, ${\mathfrak{so}}_{6,2}^2$, ${\mathfrak{so}}_{7,1}^2$ are also represented by the algebras $\varphi({\mathfrak{so}}_{5,3}^1)$, $\varphi({\mathfrak{so}}_{6,2}^1)$, $\varphi({\mathfrak{so}}_{7,1}^1)$, respectively. We note that the algebra ${\mathfrak{so}}_{4,4}$ is special in the sense that the equivalence classes of the corresponding subalgebras ${\mathfrak{so}}_{4,4}^1$ and ${\mathfrak{so}}_{4,4}^2$ coincide with that of the subalgebra ${\mathfrak{so}}_{4,4}^0$. 

Below we give explicit formulas for one such automorphism $\theta$. They can be derived from the formulas given in \cite{BB} (see pp. 655--658). Denote by $j_{rs}$, with $1\le r<s\le 8$, the $8\times 8$-matrix whose $(r,s)$th entry is $1$, the $(s,r)$th entry is $-1$, and all other entries are zero. These matrices form a basis of ${\mathfrak{so}}_8(\CC)$, and we will write $\theta$ by means of specifying the images of the basis matrices as follows
\begin{equation}
\begin{array}{ll}
\displaystyle j_{12}\mapsto (j_{12}-j_{34}+j_{56}-j_{78})/2,& \displaystyle j_{13}\mapsto (j_{13}+j_{24}+j_{57}+j_{68})/2,\\
\vspace{-0.3cm}\\
 \displaystyle j_{14}\mapsto-(j_{14}-j_{23}-j_{58}+j_{67})/2, & \displaystyle j_{15}\mapsto (j_{18}-j_{27}+j_{36}-j_{45})/2,\\
 \vspace{-0.3cm}\\
\displaystyle j_{16}\mapsto (j_{17}+j_{28}+j_{35}+j_{46})/2, & \displaystyle j_{17}\mapsto -(j_{16}+j_{25}-j_{38}-j_{47})/2,\\
\vspace{-0.3cm}\\
\displaystyle j_{18}\mapsto -(j_{15}-j_{26}-j_{37}+j_{48})/2, & \displaystyle j_{23}\mapsto -(j_{14}-j_{23}+j_{58}-j_{67})/2,\\
 \vspace{-0.3cm}\\
\displaystyle j_{24}\mapsto -(j_{13}+j_{24}-j_{57}-j_{68})/2, & \displaystyle j_{25}\mapsto (j_{17}+j_{28}-j_{35}-j_{46})/2,\\
\vspace{-0.3cm}\\
\displaystyle j_{26}\mapsto -(j_{18}-j_{27}-j_{36}+j_{45})/2, & \displaystyle j_{27}\mapsto (j_{15}-j_{26}+j_{37}-j_{48})/2,\\
\vspace{-0.3cm}\\
\displaystyle j_{28}\mapsto -(j_{16}+j_{25}+j_{38}+j_{47})/2, & \displaystyle j_{34}\mapsto (j_{12}-j_{34}-j_{56}+j_{78})/2,\\
\vspace{-0.3cm}\\\displaystyle j_{35}\mapsto -(j_{16}-j_{25}-j_{38}+j_{47})/2, & \displaystyle j_{36}\mapsto -(j_{15}+j_{26}-j_{37}-j_{48})/2,\\
\vspace{-0.3cm}\\
\displaystyle j_{37}\mapsto -(j_{18}+j_{27}+j_{36}+j_{45})/2, & \displaystyle j_{38}\mapsto -(j_{17}-j_{28}+j_{35}-j_{46})/2,\\
\vspace{-0.3cm}\\
\displaystyle j_{45}\mapsto -(j_{15}+j_{26}+j_{37}+j_{48})/2, & \displaystyle j_{46}\mapsto (j_{16}-j_{25}+j_{38}-j_{47})/2,\\
\vspace{-0.3cm}\\
\displaystyle j_{47}\mapsto (j_{17}-j_{28}-j_{35}+j_{46})/2, & \displaystyle j_{48}\mapsto -(j_{18}+j_{27}-j_{36}-j_{45})/2,\\
\vspace{-0.3cm}\\
\displaystyle j_{56}\mapsto (j_{12}+j_{34}-j_{56}-j_{78})/2, & \displaystyle j_{57}\mapsto (j_{13}-j_{24}-j_{57}+j_{68})/2,\\
\vspace{-0.3cm}\\
\displaystyle j_{58}\mapsto -(j_{14}+j_{23}-j_{58}-j_{67})/2, & \displaystyle j_{67}\mapsto (j_{14}+j_{23}+j_{58}+j_{67})/2,\\
\vspace{-0.3cm}\\
\displaystyle j_{68}\mapsto(j_{13}-j_{24}+j_{57}-j_{68})/2, & \displaystyle j_{78}\mapsto -(j_{12}+j_{34}+j_{56}+j_{78})/2.
\end{array}\label{maptheta}
\end{equation}

For $p+q=8$, $p\ge q>1$ we denote by $PSO^j_{p,q}$ the connected subgroup of $PSO_8(\CC)$ corresponding to the subalgebra ${\mathfrak{so}}_{p,q}^j$ for $j=1,2$. In order to obtain a classification of manifolds $X$ for ${\mathfrak g}_0={\mathfrak{so}}_{p,q}$, one needs to understand the orbit structure of the action of $PSO^j_{p,q}$ on ${\mathcal Q}_6$. Note that $PSO^1_{p,q}$ and $PSO^2_{p,q}$ are conjugate in $PO_8(\CC)$, and therefore the orbit structures for these two groups are equivalent.   

\subsubsection{} Assume first that ${\mathfrak g}_0={\mathfrak{so}}_{6,2}$. The main observation in this case is that the algebra ${\mathfrak{so}}_{6,2}^1$ lies in the same equivalence class as the algebra ${\mathfrak{so}}^*_8$ (recall that ${\mathfrak{so}}_{6,2}$ and ${\mathfrak{so}}^*_8$ are isomorphic), hence ${\mathfrak{so}}_{6,2}^2$ lies in the same equivalence class as $\widehat{{\mathfrak{so}}}^*_8$ (see \cite{BB}, pp. 655-658), and therefore $PSO^1_{6,2}$ and $PSO^2_{6,2}$ are conjugate in $PSO_8(\CC)$ to $PSO^*_8$ and $\widehat{PSO}^*_8$, respectively (see Example \ref{exsostar}). The classification is then given by the following theorem.    

\begin{theorem}\label{so62} Let $X$ be a connected complex manifold of dimension $6$ admitting an almost effective action by holomorphic transformations of a connected Lie group $G_0$ with Lie algebra ${\mathfrak{so}}_{6,2}$. Then $X$ is biholomorphic to one of ${\mathcal Q}_6$, ${\mathcal Q}_6\setminus{\mathcal S}_{6,2}^1$, $\Omega_{6,\,2}^{+}$, $\Omega_{6,\,2}^{-}$, $E_6^{+}$, $E_6^{-}$.
\end{theorem}

The proof is obtained by suitably combining the proofs of Theorems \ref{sopq} and \ref{sostar}.

\subsubsection{} Assume next that ${\mathfrak g}_0={\mathfrak{so}}_{7,1}$. It turns out that each of $PSO^1_{7,1}$, $PSO^2_{7,1}$ acts transitively on ${\mathcal Q}_6$ (see \cite{Wo2}, Corollary 1.7 for $r=3$, and \cite{O}, Table 11, p. 270). This fact yields the following result.   

\begin{theorem}\label{so71} Let $X$ be a connected complex manifold of dimension 6 admitting an almost effective action by holomorphic transformations of a connected Lie group $G_0$ with Lie algebra ${\mathfrak{so}}_{7,\,1}$. Then $X$ is biholomorphic to one of ${\mathcal Q}_6$, ${\mathcal Q}_6\setminus{\mathcal S}_{7,1}^1$.
\end{theorem}

\subsubsection{} Assume finally that ${\mathfrak g}_0={\mathfrak{so}}_{5,3}$. We will now determine the orbits of the $PSO_{5,3}^1$-action on ${\mathcal Q}_6$. Using formulas (\ref{maptheta}) for the triality automorphism $\theta$, one can find the algebra ${\mathfrak{so}}_{5,3}^1$ explicitly. It consists of skew-symmetric matrices $(a_{ml})_{m,l=1,\dots,8}$ with complex entries satisfying the following conditions
$$
\begin{array}{ll}
\hbox{Re}\,\left(a_{12}+a_{34}-a_{56}-a_{78}\right)=0, & \hbox{Re}\,\left(a_{14}+a_{23}-a_{58}-a_{67}\right)=0,\\
\vspace{-0.3cm}\\
\hbox{Re}\,\left(a_{13}-a_{24}-a_{57}+a_{68}\right)=0, & 
\hbox{Re}\,a_{15}=\hbox{Re}\,a_{26}=\hbox{Re}\,a_{37}=\hbox{Re}\,a_{48}, \\
\vspace{-0.3cm}\\
\hbox{Re}\,a_{18}=-\hbox{Re}\,a_{27}=\hbox{Re}\,a_{36}=-\hbox{Re}\,a_{45}, &
\hbox{Re}\,a_{16}=-\hbox{Re}\,a_{25}=-\hbox{Re}\,a_{38}=\hbox{Re}\,a_{47},\\
\vspace{-0.3cm}\\
\hbox{Re}\,a_{17}=\hbox{Re}\,a_{28}=-\hbox{Re}\,a_{35}=-\hbox{Re}\,a_{46}, &
\hbox{Im}\,\left(a_{17}+a_{28}-a_{35}-a_{46}\right)=0, \\ 
\vspace{-0.3cm}\\
\hbox{Im}\,\left(a_{18}-a_{27}+a_{36}-a_{45}\right)=0, &
\hbox{Im}\,\left(a_{16}-a_{25}-a_{38}+a_{47}\right)=0, \\
\vspace{-0.3cm}\\
\hbox{Im}\,\left(a_{15}+a_{26}+a_{37}+a_{48}\right)=0.
\end{array}
$$
A straightforward (but lengthy) calculation using the above conditions yields that there are exactly two $PSO_{5,3}^1$-orbits in ${\mathcal Q}_6$: a 9-dimensional closed orbit $\Gamma$ and the open orbit ${\mathcal Q}_6\setminus\Gamma$. The orbit $\Gamma$ contains, for example, all points that are represented in homogeneous coordinates in $\PP^7$ by vectors $(c,id)$, where $c,d\in\RR^4$, $||c||=||d||=1$, and all points represented by vectors $(c,i\overline{c})$, where $c=(c_1,c_2,c_3,c_4)\in\CC^4$ is such that $||c||=1$ and $\hbox{Im}\left(\sum_{j=1}^4c_j^2\right)=0$.

Arguing as in the proofs of Theorems \ref{slr}, \ref{sopq}, we now obtain the following result.

\begin{theorem}\label{so53} Let $X$ be a connected complex manifold of dimension 6 admitting an almost effective action by holomorphic transformations of a connected Lie group $G_0$ with Lie algebra ${\mathfrak{so}}_{5,\,3}$. Then $X$ is biholomorphic to one of ${\mathcal Q}_6$, ${\mathcal Q}_6\setminus{\mathcal S}_{5,3}^1$, $\Omega_{5,\,3}^{+}$, $\Omega_{5,\,3}^{-}$, ${\mathcal Q}_6\setminus\Gamma$.
\end{theorem}

Theorems \ref{sopq}, \ref{sopqn=3}, \ref{sostar}, \ref{so62}, \ref{so71}, \ref{so53} complete the proof of Theorem \ref{maintheoremI}.

Let us conclude our discussion of the case ${\mathfrak g}_0={\mathfrak{so}}_{5,3}$ by sketching another approach to understanding actions of $PSO^j_{p,q}$ on ${\mathcal Q}_6$ for $j=1,2$. This replaces considerations of triality with the study of actions of $PSO^{\circ}_{p,q}$ on Grassmannians associated to ${\mathcal Q}_6$.  

The Dynkin diagram of ${\mathfrak{so}}_8(\CC)$ has three outer nodes which are permuted by triality automorphisms. Each of these nodes corresponds to a $PSO_8(\CC)$-homogeneous manifold equivalent to ${\mathcal Q}_6$ endowed with one of three actions of $PSO_8(\CC)$. If the left-most node corresponds to the standard action, then each of the other two nodes corresponds to an action defined by composing the standard one with a triality automorphism of $PSO_8(\CC)$. Each of the two non-standard actions on ${\mathcal Q}_6$ is equivalent to the standard action of $PSO_8(\CC)$ on a family of linearly embedded 3-dimensional projective subspaces in ${\mathcal Q}_6$. If we view ${\mathcal Q}_6$ as the space of $[\cdot ,\cdot]_{p,q}$-isotropic 1-dimensional subspaces of $\CC^8$ (see (\ref{zzpq}) for the definition of $[\cdot ,\cdot]_{p,q}$), then these are the two $SO_8(\CC)$-homogeneous families of 4-dimensional $[\cdot ,\cdot]_{p,q}$-isotropic subspaces. The action of each group $PSO_{p,q}^j$ on ${\mathcal Q}_6$ has the same orbit structure as the associated action of $SO_{p,q}^{\circ}$ on one of these families. Let us describe these associated actions in concrete terms in the case $p=5$, $q=3$ leaving the other cases to the reader.

Let $W$ be a given $[\cdot ,\cdot]_{5,3}$-isotropic 4-plane and note that it contains a 1-dimensional subspace $L$ which is
positive with respect to the Hermitian form $\langle \cdot ,\cdot\rangle_{5,3}$ (see (\ref{anglezzpq})). By applying an element of $SO_{5,3}^{\circ}$ we may assume that $L$ is generated by $v:=(1,i,0,\ldots ,0)$. Now write $W=L\oplus P$ where $P$ is the $\langle \cdot ,\cdot\rangle_{5,3}$-orthogonal complement of $L$. Observe that since $W$ is $[\cdot ,\cdot]_{5,3}$-isotropic, it follows that the 3-plane $P$ is $\langle \cdot ,\cdot\rangle_{5,3}$-orthogonal to
$E:=\hbox{Span}\,\{v,\overline{v}\}$. Note also that the stabilizer of $E$ acts as $SO_{3,3}^{\circ}$ on the orthogonal complement of $E$ in $\CC^8$. Thus we must only describe the action of $SO_{3,3}^{\circ}$ on families of 3-planes in $\CC^6$ which are isotropic with respect to $[\cdot,\cdot]_{3,3}$. 

Examples of two such families are those of the 3-planes which are isotropic with respect to $\langle \cdot ,\cdot\rangle_{3,3}$. These define the closed $SO_{5,3}^{\circ}$-orbits in each of the families of $[\cdot ,\cdot]_{5,3}$-isotropic 4-planes in $\CC^8$. If $(a,b,c)$ denotes the Hermitian signature of a 4-plane $W$, with $a$ the dimension of a maximal positive subspace, $b$ the dimension of a maximal negative subspace and $c$ the dimension of the degeneracy, then we see that this closed orbit is the manifold of 4-planes of signature $(1,0,3)$. An explicit computation at a base point provided by the above description shows that it is 9-dimensional.

If the 3-plane $P$ arising in the decomposition $W=L\oplus P$ is not $\langle \cdot ,\cdot\rangle_{3,3}$-isotropic, then the above-described procedure of splitting off standard 2-planes $E$ shows that $P$ is of signature $(1,1,1)$ and can be moved by $SO_{3,3}^{\circ}$ to a space spanned by vectors of the form $v_+=(1,\pm i,0,0,0,0)$, $v_0=(0,0,1,\pm 1,0,0)$, $v_-=(0,0,0,0,1,\pm i)$. In general, the signs may play a role. For example, this is the case for the $SO_{4,4}^{\circ}$-action where there are in fact two orbits of signature $(2,2,0)$ in each family. However, in the present situation we can apply simple involutions in the maximal compact group $SO_3(\RR)\times SO_3(\RR)$, which is naturally embedded in $SO_{3,3}^{\circ}$, to bring the 3-tuple $(v_+,v_0,v_{-})$ to one of two normal forms. These correspond to orbits of signature $(2,1,1)$ in the two families of isotropic 4-planes. Thus in each family the group $SO_{5,3}^{\circ}$ has two orbits, an open orbit of signature $(2,1,1)$ and a closed orbit of signature $(1,0,3)$. The closed orbit is the one denoted
by $\Gamma$ above.

Fakult\"at f\"ur Mathematik, Ruhr-Universit\"at Bochum, Universit\"atsstra\ss e 150, D-44801 Bochum, Germany \hspace{0.1cm}$\bullet$\hspace{0.1cm} {\tt e-mail: ahuck@cplx.ruhr-uni-bochum.de}
\vspace{1cm}

Department of Mathematics, The Australian National University, Canberra, ACT 0200, Australia \hspace{0.1cm}$\bullet$\hspace{0.1cm} {\tt e-mail: alexander.isaev@maths.anu.edu.au}

\end {document}